\documentclass[11pt]{amsart}
 \usepackage{etex}
 \usepackage{relsize}
 \usepackage{amsmath, amscd, amsthm, amssymb, graphics, mathrsfs, setspace,fancyhdr,geometry,tabularx,shapepar, hyperref, xcolor, tikz, booktabs, dsfont,amsfonts, marvosym, amsbsy, bm, tikz, array,}
 \usetikzlibrary{matrix,arrows,backgrounds,cd}
 \usepackage{mathbbol}
 \usepackage{youngtab}
 \usepackage{young}
 \usepackage{ytableau}
 \usetikzlibrary{decorations.pathreplacing}
 \newcommand{\tikznode}[3][inner sep=0pt]{\tikz[remember
 	picture,baseline=(#2.base)]{\node(#2)[#1]{$#3$};}}
 
 \usepackage{lineno}  
 \modulolinenumbers[5]

 \makeatletter
 \newcommand{\adjunction}{\@ifstar\named@adjunction\normal@adjunction}
 \newcommand{\normal@adjunction}[4]{%
 	#1\colon #2%
 	\mathrel{\vcenter{%
 			\offinterlineskip\m@th
 			\ialign{%
 				\hfil$##$\hfil\cr
 				\longrightharpoonup\cr
 				\noalign{\kern-.3ex}
 				\smallbot\cr
 				\longleftharpoondown\cr
 			}%
 	}}%
 	#3 \noloc #4%
 }
 \newcommand{\named@adjunction}[4]{%
 	#2%
 	\mathrel{\vcenter{%
 			\offinterlineskip\m@th
 			\ialign{%
 				\hfil$##$\hfil\cr
 				\scriptstyle#1\cr
 				\noalign{\kern.1ex}
 				\longrightharpoonup\cr
 				\noalign{\kern-.3ex}
 				\smallbot\cr
 				\longleftharpoondown\cr
 				\scriptstyle#4\cr
 			}%
 	}}%
 	#3%
 }
 \newcommand{\longrightharpoonup}{\relbar\joinrel\rightharpoonup}
 \newcommand{\longleftharpoondown}{\leftharpoondown\joinrel\relbar}
 \newcommand\noloc{%
 	\nobreak
 	\mspace{6mu plus 1mu}
 	{:}
 	\nonscript\mkern-\thinmuskip
 	\mathpunct{}
 	\mspace{2mu}
 }
 \newcommand{\smallbot}{%
 	\begingroup\setlength\unitlength{.15em}%
 	\begin{picture}(1,1)
 	\roundcap
 	\polyline(0,0)(1,0)
 	\polyline(0.5,0)(0.5,1)
 	\end{picture}%
 	\endgroup
 }
 \makeatother
 \hypersetup{colorlinks=false}
 \theoremstyle{plain}
 \newtheorem{thm}{Theorem}[section]
 \newtheorem{lem}[thm]{Lemma}
 
 \newtheorem{cor}[thm]{Corollary}
 \newtheorem{prop}[thm]{Proposition}

 \newtheorem*{thm*}{Theorem}
 \newtheorem*{problem*}{Problem}

 \theoremstyle{definition}
 \newtheorem{defn}[thm]{Definition}
 
 \newtheorem{notation}[thm]{Notation}
 
 \newtheorem{rem}[thm]{Remark}

 \newtheorem{ex}[thm]{Example}

 \newtheorem*{introconjecture}{Conjecture}
  \newtheorem*{introremark}{Remark}
 \newtheorem{introtheorem}{Theorem}

 \numberwithin{equation}{section}

 \newcommand{\Spec}{\operatorname{Spec}}

 \renewcommand{\O}{{\mathcal O}}
 \renewcommand{\hom}{\operatorname{Hom}}
 \newcommand{\kk}{\mathbb{k}}
 
 \newcommand{\cplx}{{\mathbb C}}

 \newcommand{\Nat}{\mathbb{N}}

 \newcommand{\im}{\operatorname{Im}}
 
 \newcommand{\ra}{\rightarrow}

 \newcommand{\gm}{\mathbb{G}_{m}}

 \newcommand{\op}[1]{\operatorname{#1}}
 \newcommand{\spec}{\op{Spec}}
 \newcommand{\coker}{\operatorname{coker}}

 \newcommand*{\sheafhom}{\mathscr{H}\kern -.5pt om}
 
 \newcommand{\Omodmod}{/ \! \! / 0}

 \newcommand{\kapreezer}{\raisebox{.448em}{-}{\kern-.246em\raisebox{.3em}{-}}{\kern-0.2399em\mathsf{K}}{\kern-0.165em \raisebox{.449em}{-}}{\kern-0.439em \raisebox{.3em}{-}}}
 \newcommand{\lweezer}{\raisebox{.44em}{-}{\kern-.29em\raisebox{.3em}{-}}{\kern-0.18em\mathsf{K}}{\kern-0.24em \raisebox{.449em}{-}}{\kern-0.439em \raisebox{.3em}{-}}}

 \newcommand{\rox}{Z\Omodmod}
 \newcommand{\ktimes}{\times_{\kk}}

 \newcommand{\gl}{G}
 \newcommand{\inv}{^{-1}}

 \newcommand{\sch}{\mathcal{O}}
 \newcommand{\sQ}{\mathbf{Q}}
 
 \newcommand{\lan}{(}
 \newcommand{\ran}{)}

 \newcommand{\hpi}{\pi^\sharp}
 \newcommand{\hsi}{\sigma^\sharp}
 \newcommand{\kg}{\kk\left[\gl\right]}

 \newcommand{\pQ}{\widehat{Q}}
 \newcommand{\pQs}{\pQ^{+}}

 \newcommand{\mdx}{\widetilde{\Delta}_Z}
 \newcommand{\qsquare}{Q_Z\hspace{1mm}_s\hspace{-1mm}\otimes_p\hspace{.5mm}Q_Z}
 \newcommand{\kten}{\otimes_\kk}
 
 \newcommand{\bZ}{Z}
 \newcommand{\Db}{\mathsf{D^b}}

 \DeclareMathOperator{\Hom}{Hom}
 
 \newcommand{\MB}[1]{}
 \newcommand{\PM}[1]{}
 \newcommand{\RV}[1]{} 
 \newcommand{\DF}[1]{}
 \newcommand{\NC}[1]{}
 
 \numberwithin{equation}{section}
 \numberwithin{thm}{subsection}
 
\begin{document}

 \title{Kernels for Grassmann flops}

 
 	\author[Ballard]{Matthew R. Ballard}
 \address{
 	\begin{tabular}{l}
 		Matthew Robert Ballard \\
 		\hspace{.1in} University of South Carolina\\
 		\hspace{.1in} Department of Mathematics \\
 		\hspace{.1in} Columbia, SC, USA  \\
 		\hspace{.1in} Email: {\bf ballard@math.sc.edu} \\
 		\hspace{.1in} Webpage: \url{https://www.matthewrobertballard.com/} \\
 	\end{tabular}
 }

 \author[Chidambaram]{Nitin K. Chidambaram}
 \address{
 	\begin{tabular}{l}
 		Nitin K. Chidambaram \\
 		\hspace{.1in} Max Planck Institut f\"ur Mathematik, \\
 		\hspace{.1in} Bonn, Germany\\
 		\hspace{.1in} Email: {\bf kcnitin@mpim-bonn.mpg.de} \\
 		\hspace{.1in} Webpage: \url{https://guests.mpim-bonn.mpg.de/kcnitin} \\
 	\end{tabular}
 }

 \author[Favero]{David Favero}
 \address{
 	\begin{tabular}{l}
 		David Favero \\
 		\hspace{.1in} University of Alberta \\
 		\hspace{.1in} Department of Mathematical and Statistical Sciences \\
 		\hspace{.1in}  Edmonton, AB, Canada  \\
 		\\
 		\hspace{.1in} Korea Institute for Advanced Study \\
 		\hspace{.1in}   Seoul, Republic of Korea  \\
 		\hspace{.1in} Email: {\bf favero@ualberta.ca} \\
 		\hspace{.1in} Webpage: \url{https://sites.ualberta.ca/~favero/} \\
 	\end{tabular}
 }

 \author[McFaddin]{Patrick K. McFaddin}
 \address{
 	\begin{tabular}{l}
 		Patrick K. McFaddin \\ 
 		\hspace{.1in} Fordham University\\
 		\hspace{.1in} Department of Mathematics \\
 		\hspace{.1in} New York, NY, USA  \\
 		\hspace{.1in} Email: {\bf pmcfaddin@fordham.edu} \\
 		\hspace{.1in} Webpage: \url{http://mcfaddin.github.io/} \\
 	\end{tabular}
 }
 
 \author[Vandermolen]{Robert R. Vandermolen}
 \address{
 	\begin{tabular}{l}
 		Robert Richard Vandermolen \\
 		\hspace{.1in} Saint Mary-of-the-Woods College\\
 		\hspace{.1in}  St. Mary-of-the-Woods, IN, USA \\
 		\hspace{.1in} Email: {\bf Robert.Vandermolen@smwc.edu} \\
 		\hspace{.1in} Webpage: \url{https://robertvandermolen.github.io/} \\
 	\end{tabular}
 }

\begin{abstract}
 We develop a generalization of the $Q$-construction of the first author, Diemer, and the third author for Grassmann flips. This generalization provides a canonical idempotent kernel on the derived category of the associated global quotient stack. The  idempotent kernel, after restriction, induces a semi-orthogonal decomposition which compares the flipped varieties. Furthermore its image, after restriction to the geometric invariant theory semistable locus, ``opens'' a canonical ``window'' in the derived category of the quotient stack. We check this window coincides with the set of representations used by Kapranov to form a full exceptional collection on Grassmannians. 
\end{abstract}

\keywords{ Derived categories, birational geometry, representation theory}

\ytableausetup{centertableaux}
\maketitle
\addtocounter{section}{1}




\section*{Introduction}

Derived categories, once viewed as a mere technical book-keeping device, have flourished as a topic of investigation as volumes of literature have exposed their geometric nature. Derived categories of coherent sheaves on algebraic varieties bind algebraic geometry to commutative algebra, representation theory, symplectic geometry, and theoretical physics in deep and surprising ways. 

These bindings come in the form of fully-faithful functors or, better yet, equivalences relating different varieties or categories. An obvious and central question: what is a reasonably robust and general source for such functors?
Experience in algebraic geometry tells us that moduli spaces are often a good place to look but beyond this source the examples of fully-faithful functors are more idiosyncratic. 

Recently, a new construction, in the context of group actions, was introduced in \cite{BDF}, which we call the $Q$-construction. Given a variety $X$ with a $\mathbb{G}_m$-action, the authors constructed an idempotent kernel on the equivariant derived category $\Db([X/\mathbb{G}_m])$. The kernel $Q$, being the identity on its essential image, fully-faithfully identifies an interesting component of the derived category $\Db([X/\mathbb{G}_m])$. In fact, it always gives a two-term semi-orthogonal decomposition. This construction has some natural extensions. 

Following Drinfeld \cite{drin}, we can recognize it as a piece of a more general story. The inclusion $\mathbb{G}_m \subset \mathbb{A}^1$ can be viewed as a partial compactification of $\mathbb{G}_m$ as a monoid in schemes. The fibers of the multiplication map $\mathbb{A}^1 \times \mathbb{A}^1 \to \mathbb{A}^1$ are a family of $\mathbb{G}_m$ orbits which degenerate over $0 \in \mathbb{A}^1$. From Drinfeld's perspective, the idempotent kernel constructed in \cite{BDF} is the structure sheaf of a variety that parametrizes such degenerations in $X$. 

This viewpoint allows for an immediate generalization: we can replace $\mathbb{A}^1$ with $M$ where $M$ is monoidal scheme. If we have a variety with an action of the units of $M$, we can produce a kernel for $X$. In this paper, we study the monoidal scheme $\op{End}(V)$ for $V$ a finite dimensional vector space and the natural action of the units $\operatorname{GL}(V)$ on the vector space $Z = \op{Hom}(V,W) \times \op{Hom}(W',V)$ where $W,W'$ are two additional finite dimensional vector spaces. 

In this setting, we construct an object $Q \in \Db([Z \times Z/\op{GL}(V) \times \op{GL}(V)])$ and show that the idempotent property still holds:

\begin{introtheorem}
 There exists a morphism of kernels $Q \to \Delta$ inducing an isomorphism of $Q \circ Q \to Q$, where $\circ$ denotes convolution of kernels and $\Delta$ is the kernel of the identity.  
\end{introtheorem} 

The robustness of the equivariant setting is indicated by the fact that any flip between normal varieties arises as a variation of GIT problem for a $\mathbb{G}_m$ action on a variety. We remind the reader of the following conjecture of Bondal and Orlov \cite{BO} extended by Kawamata \cite{Kaw}:

\begin{introconjecture}[Bondal-Orlov 1995]
 Assume that $X$ and $X^\prime$ are smooth complex varieties. If $X$ and $X^\prime$ are related by a flop, then there is a $\cplx$-linear triangulated equivalence of their bounded derived categories of coherent sheaves
 \begin{displaymath}
  \Db(X) \cong \Db(X^\prime).
 \end{displaymath}
\end{introconjecture}

As an application of the $Q$-construction, in \cite{BDF}, the first author, Diemer, and the third author gave a means of constructing a kernel on $X \times X^\prime$ for any $D$-flip $(X,X^\prime)$. For flops of smooth projective varieties, the associated integral transform is conjectured to be the desired equivalence of Bondal and Orlov. This provides a single unified, though still conjectural, approach to constructing equivalences from flops.

In the setting of this paper, variation of GIT between $Z^+,Z^-$, the two GIT quotients of $Z$,  amounts to a flip which was studied by Donovan and Segal \cite{DS} when $W$ and $W'$ have the same dimension and $\kk$ is an algebraically-closed base field of characteristic zero.  This is called a Grassmann flop since it comes from contracting the zero section of a vector bundle over a Grassmannian.  Donovan and Segal exhibited equivalences for Grassmann flops using a set of representations identified by Kapranov \cite{Kapp}. In \cite{BLV}, Buchweitz, Leuschke, and Van den Bergh showed that structure sheaf of fiber product $\mathcal O_{Z^+ \times_Z Z^-}$ for Grassmann flop is a kernel of the equivalence and they showed that, for the analogous flip setting, it provides admissible embeddings of derived categories. 

As an application to these conjectures, we show the $Q$-construction descends to provide an appropriate kernel for Grassmann flips. 
\begin{introtheorem} \label{thm:intro_kernel}
 Let $\kk$ be an (arbitrary) field of characteristic zero. 
 For a Grassmann flip, the $Q$ construction, restricted to the semi-stable loci, 
 induces the following semi-orthogonal decompositions:
 \begin{enumerate}
 \item If $\dim W > \dim W'$, then there is a semi-orthogonal decomposition
 \[
\Db(Z^+)  = \langle \mathsf O_{ \dim V,\dim W-1}, \Db(Z^-) \rangle
 \]
 where $\mathsf O_{\dim V,\dim W-1}$ is a category of explicitly described objects supported on the unstable locus for $Z^-$ (see in Definition~\ref{def:O}).
  \item If $\dim W < \dim W'$, then there is a semi-orthogonal decomposition
 \[
\Db(Z^-)  = \langle \mathsf O_{ \dim V,\dim W'-1}, \Db(Z^+) \rangle
 \]
  where $\mathsf O_{\dim V,\dim W'-1}$ is a category of explicitly described objects supported on the unstable locus for $Z^+$ (see Definition~\ref{def:O}).
  \item If $\dim W = \dim W'$, then there is an equivalence
 \[
\Db(Z^+)  \cong \Db(Z^-).
 \]
 \end{enumerate}
\end{introtheorem}

We also check, a posteieri, that the kernel from the $Q$-construction agrees 
with the fiber product. Additionally, we show that the equivalence holds for any 
twisted Grassmann flop over a general field in characteristic $0$.

Finally, the utility of ``windows'' in equivariant derived categories is evident from works such as, for example \cite{HL,BFK,DS,SV}. Identifying windows involves some choices or special conditions. The next result indicates that windows come simply from the choice of the monoid $M$ compactifying $G$. 

\begin{introtheorem} \label{thm:intro_window}
 Let $Q_+ := Q|_{Z^+ \times Z}$. Then 
 \begin{itemize}
  \item The functor 
 \begin{displaymath}
  \Phi_{Q_+} : \Db(Z^+) \to \Db\left( \left[Z / \operatorname{GL}(V) \right] \right)
 \end{displaymath}
 is fully-faithful.
  \item The restriction map $j^\ast : \Db\left( \left[Z / \operatorname{GL}(V) \right] \right) \to \Db(Z^+)$ is a left inverse to $\Phi_{Q_+}$.
  \item Kapranov's representations form a set of generators for the essential image of $\Phi_{Q_+}$. 
 \end{itemize}
\end{introtheorem}

Theorem~\ref{thm:intro_window}, in particular, provides a completely geometric explanation for the appearance of Kapranov's representations.  Our method of monoid compactification can therefore be seen as part of a program to produce canonical windows for quotients via linearly reductive groups.  

Theorems~\ref{thm:intro_kernel} and~\ref{thm:intro_window} demonstrate that, 
in the case of Grassmann flips, the $Q$-construction plays a unifying role in 
understanding the relationship between birational geometry and derived categories. 

\begin{introremark}
 The authors became aware of \cite{BLV} after the submission of this article. We thank Michel Van den Bergh for pointing it out. 
\end{introremark}

\section{Notation and Conventions}\label{Background}
\noindent

 Throughout, $\kk$ denotes a field of characteristic zero.  
For $\ell\in\Nat$, $\ell\neq0$, we let $[\ell]:=\{1,...,\ell\}$.
 We denote by $\mathsf{Vec}_{\kk}$ the category of finite-dimensional $\kk$-vector spaces.
 We utilize standard results in Geometric Invariant Theory and use the notation of \cite{GIT} as much as possible. All schemes considered here are $\kk$-schemes.  For a $\kk$-scheme $Z$, we denote its ring of regular functions by $\kk[Z]$. The word \emph{point} will always mean a $\kk$-\emph{point}.

Throughout, fix a $d$-dimensional $\kk$-vector space  $V$, a $ m $-dimensional vector space  $W$ and a $ m' $-dimensional vector space such that $\dim W, \dim W' \geq \dim V$. 
For any $\kk$-vector space $U$ we can associate a scheme over $\spec(\kk)$ defined as the spectrum of the symmetric algebra of the dual space $U^\vee$, that is  
\begin{equation*}
\spec\left( \op{Sym} \left(U^\vee\right)\right).
\end{equation*}
More generally, for a $\kk$-scheme $X$ and a locally free $\sch_X$-module $M$ we can consider the relative spectrum of the symmetric sheaf of algebras of $M^\vee$, which we refer to as the {\bf total space} of $M$.

We fix the group $G$ to be the general linear algebraic group
\begin{equation*}
G = \op{GL}(V) :=\spec\left(\kk\left[C,\left(\det\left(C\right)\right)\inv \right] \right),
\end{equation*}
where $C=\left(c_{ij}\right)_{i,j\in[d]}$ is a collection of indeterminates. We use $\det\left(C\right)$ to denote the polynomial
\begin{equation*}
\det\left( C\right):=\sum_{\delta\in S_n}(-1)^{\op{sgn}(\delta)}\left(\prod_{i\in[d]}c_{i\delta(i)}\right).
\end{equation*} 

Next, recall that for a $\kk$-scheme $Z$, an action of $\gl$ on $Z$ is defined by a morphism of schemes 
\begin{equation*}
\sigma_Z:\gl\ktimes Z\ra Z.
\end{equation*}
If $Z$ and $Y$ are $\kk$-schemes with $\gl$-actions given by $\sigma_Z$ and $\sigma_Y$, we say that a morphism $f:Z\ra Y$ is $\gl$-equivariant whenever the following diagram commutes:
\begin{center}\begin{tikzcd}
\gl\ktimes Z\ar[rrr,"1_{\gl}\ktimes f"]\ar[dd,"\sigma_{Z}"']&&&\gl\ktimes Y\ar[dd,"\sigma_{Y}"]\\
&&&&\\
Z\ar[rrr,"f"]&&&Y
\end{tikzcd}\end{center}
Most of this work deals with categories whose objects  carry a $\gl$-action and whose morphisms are $\gl$-equivariant. To denote such categories, we simply use the superscript $\gl$. For example, let $R$ be a commutative ring with a $\gl$-action. Then, we denote the category of $\gl$-equivariant modules by $\op{Mod}^{\gl}(R)$.

In practice, we will consider the case where $R$ is a polynomial ring and think of the coordinates along with the $\gl$-action as follows. Given a collection $A=\{a_{ij}\}_{i\in[m],j\in[d]}$ of indeterminates, for any subsets $J\subseteq[d]$, $I\subseteq[m]$, we let $A_{I, J}$ denote the collection of variables 
\begin{equation*}
A_{I,J}:=(a_{ij})_{i\in I,\; j\in J}.
\end{equation*}
For $I\subseteq[m]$ with $|I|=d$, we may list this set in increasing order, and denote the corresponding ordered set by $I:=\{\ell_1,..,\ell_d\}$. We then write 
\begin{equation*}
\det\left(A_{I,[d]}\right):=\sum_{\sigma\in S_d}\op{sgn}(\sigma)\left(\prod_{1\leq i\leq d}A_{\ell_i,\sigma(i)}\right).
\end{equation*}
Given two collections $A=\{a_{ij}\}_{i\in[m],j\in[d]}$ and $B=\{b_{ij}\}_{i\in[d],j\in[m]}$, we use $BA$ to denote the collection of polynomials
\begin{equation*}
\left\{ \sum_{\ell=1}^mb_{i\ell}a_{\ell j}\right\}_{i\in[d], j\in[d]}.
\end{equation*}
Lastly, given two $\kk$-algebras $R$ and $S$ and elements $r \in R$ and $s \in S$, we let $r^L$ and $s^R$ denote the elements $r \otimes 1$ and $1 \otimes s$ in $R \kten S$.

Throughout the text, we will use $ Z $ to denote the following affine space 
\begin{equation*}
Z := \hom_{\mathsf{Vec}_\kk}(V,W) \times \hom_{\mathsf{Vec}_\kk}(W',V),
\end{equation*} and we will denote its coordinate ring by 
\[
R:= k[Z].
\]
\section{The kernel}\label{section:q}

Recall that in \cite{BDF} the authors exhibit a kernel using a partial compactification of a certain $\gm$-action. We follow a similar line of reasoning in the case of $ \gl $-actions on spaces of a certain type.

\subsection{The Grassmann flip}

The vector space
\begin{equation*}
\hom_{\mathsf{Vec}_\kk}(V,W) \oplus \hom_{\mathsf{Vec}_\kk}(W',V),
\end{equation*}
 carries a natural action of $ \gl $; for  $\varsigma\in\gl$ (a point), $\varphi\in\hom_{\mathsf{Vec}_{\kk}}(V,W)$, and $\vartheta\in\hom_{\mathsf{Vec}_{\kk}}(W',V)$ the action is given as 
 \[
\varsigma\cdot(\varphi,\vartheta) = (\varphi\circ\varsigma,\varsigma\inv\circ\vartheta)
\]
Let us provide a few more specifics concerning the above action of $\gl$. For such an object $Z$, the induced action
\begin{equation*}
\sigma_Z:\gl\ktimes Z\ra Z
\end{equation*}
is equivalent to the co-action as Hopf algebra modules as follows.  Choosing bases for $V$, $W$ and $W'$, we may write
\begin{equation*}
Z=\spec\left(\kk\left[ \{a_{ij}\}_{i\in[d],j\in[m']},\{b_{ij}\}_{i\in[m],j\in[d]}\right] \right),
\end{equation*}
where we recall that $m=\dim W$ and $m'=\dim W'$. Letting $B: =\left(b_{ij}\right)$ and $A:=\left( a_{ij}\right)$, we have 
\begin{equation}\label{R}
Z= \spec\left(\kk\left[ A,B\right] \right).
\end{equation}
The co-action on the global sections 
\begin{equation*}
\hsi_Z:\kk[Z]\ra\kk\left[ \gl\right]\otimes_\kk \kk[Z],
\end{equation*}
is defined on the generators as
 \begin{align*}
  b_{ij}&\mapsto \left(\det(C)\right)\inv\sum_{r=1}^d \op{Adj}(C)_{rj}\otimes_\kk b_{ir}\\
  a_{ij}&\mapsto \sum_{r=1}^d c_{ir}\otimes_\kk a_{rj},
\end{align*}
where $\op{Adj}(C)$ is the adjugate matrix of $C$.

Furthermore, the projection
\begin{equation*}
\pi_Z:\gl\ktimes Z\ra Z
\end{equation*}
induces the map
\begin{equation*}
\hpi_Z:\kk[Z]\ra \kk\left[ \gl\right]\otimes_\kk\kk[Z].
\end{equation*}

In order to get our GIT problem of interest, we consider the two open sets 
\begin{align*}
U^+&:=\Big(\hom(V,W)\setminus\{\varphi:\op{rank}(\varphi)\leq(d-1)\}\Big)\oplus\hom(W',V)\\
U^-&:=\hom(V,W)\oplus\Big(\hom(W',V)\setminus\{\vartheta:\op{rank}(\vartheta)\leq (d-1)\}\Big),
\end{align*} and define the associated GIT quotients as 
\[
Z^+ := \left[ U^+/\gl\right], \qquad Z^- := \left[ U^-/\gl\right].
\] 
The resulting birational transformation between $ Z^+ $ and $ Z^- $ is called the Grassmann flip.  It is called a Grassmann flop when $ \dim W' = \dim W $.

\subsection{The object \texorpdfstring{$ \sQ $}{Q}}
Before turning attention to our object $\sQ$, we introduce the object $\Delta$, which gives the kernel of the identity functor in the equivariant setting.
\begin{notation}
Given $ Z = \Hom(V,W) \times \Hom(W',V)  $ as before, we define the following scheme
\begin{equation*}
\Delta_Z:=Z\ktimes\gl.
\end{equation*}
\end{notation}

 The scheme $ \Delta_Z $ is equipped with a natural $ \gl \ktimes \gl $-action as follows.

\begin{lem}\label{lemma:dubact}
	The scheme $\Delta_Z$ has a natural $(\gl\ktimes\gl)$-action 
	\begin{equation*}
	\sigma_{\Delta_Z}:\gl\ktimes\gl{}\ktimes\Delta_Z\ra\Delta_Z
	\end{equation*}
	uniquely determined by the co-action
	\begin{equation*}
	\hsi_{\Delta_Z}:\kk[\Delta_Z] = \kk[\gl] \otimes R \ra\kg\otimes_k\kg\otimes_\kk \kk[\Delta_Z]
	\end{equation*}
	defined by
	\begin{align*}
	\hsi_{\Delta_Z}\left(1 \otimes r \right)&= \Big(\iota_1\otimes 1_{\Delta_Z}\Big) \circ \hsi_{Z}(r)\\
	\hsi_{\Delta_Z}\left(t \otimes 1 \right)&=\left(\left(1\otimes \mu^\sharp \right) \circ \Big(\beta^\sharp\otimes 1\Big)\circ s^\sharp \circ \mu^\sharp (t)\right) \otimes 1_R.
	\end{align*}
	Here $r\in R$, $t\in\kk[\gl]$ and $\iota_1:\kg\ra\kg\otimes_\kk\kg$ is the natural inclusion into the first component;
	\begin{itemize}
		\item $\beta^\sharp:\kg\ra\kg$ is the co-inverse,
		\item $ \mu^\sharp: \kg\ra\kg \otimes \kg $ is the group co-multiplication, and
		\item  $ s^\sharp:\kg \otimes \kg\ra\kg \otimes \kg $ switches the factors in the tensor product. 
	\end{itemize}
	Moreover, the map $\pi_Z\ktimes\sigma_Z:\gl\ktimes Z\ra Z\ktimes Z$ is equivariant with respect to this $ \gl\ktimes\gl $ action.
\end{lem}
\begin{proof}
The proof is a straight-forward diagram chase and is left to the reader.
\end{proof}

Then, by using the $(\gl\ktimes\gl)$-equivariant morphism
\begin{equation*}
\pi_Z\ktimes\sigma_Z : Z\ktimes\gl \ra Z\ktimes Z,
\end{equation*}
 we get a $(\gl\ktimes\gl)$-equivariant sheaf  of modules  over $Z\ktimes Z$ associated to $\Delta_Z$, which we denote
\begin{equation*}
\mdx:=(\pi_Z\ktimes\sigma_Z)_*\sch_{\Delta_Z},
\end{equation*}
where $\sch_{\Delta_Z}$ denotes the structure sheaf of the affine scheme $\Delta_Z$.

We view $ \mdx $ as an object of $\Db(\text{Qcoh}^{\gl \ktimes \gl} Z\ktimes Z)$, and claim that it is Fourier-Mukai kernel for the identity functor on the bounded $\gl$-equivariant derived category $\Db(\text{Qcoh}^{\gl} Z)$.
\begin{lem}\label{lem:invar}
Let $Z=\spec(R)$ be before and $M$ be an object of $\op{Mod}^{\gl}(R)$. Then there is a $\kk[\gl]$-co-module isomorphism 
\begin{equation*}
\Big(\kk[\gl]\otimes_\kk M\Big)^{\gl}\cong M,
\end{equation*}
where $\kk[\gl]\kten M$ is given the left $\kk[\gl]$-co-action as a $\kk[\Delta_Z]$-module.
\end{lem}
\begin{proof}
Note that there is a natural morphism 
\begin{displaymath}
 M \to \Big(\kk[\gl]\otimes_\kk M\Big)^{\gl}
\end{displaymath}
given by the equivariant structure of $M$. Since the extension $\overline{\kk}/\kk$ is faithfully-flat, it suffices to show that this map is an isomorphism over $\overline{\kk}$. Assume that $\overline{\kk} = \kk$. 

By the Peter-Weyl Theorem, there is a decomposition $\kk[\gl]=\bigoplus S_i\otimes S_i^{\vee}$, where $S_i$ runs over every irreducible representation of $\gl$. Furthermore, since $\gl$ is linearly reductive, we have a decomposition $M=\bigoplus M_i$ into irreducible components.
Thus, we have
\begin{equation*}
\Big(\kk[\gl]\otimes_\kk M\Big)^{\gl}\cong \bigoplus S_i\otimes(S_i^{\vee}\otimes M_i)^{\gl}\cong \bigoplus S_i\otimes\hom^{\gl}_{\kk}(S_i,M_i),
\end{equation*}
and our result follows from Schur's Lemma.
\end{proof}

\begin{lem}\label{kernelid}
The object 
\begin{equation*}
\widetilde{\Delta}_Z\in \Db(\op{Qcoh}^{\gl\ktimes\gl}(Z \ktimes Z))
\end{equation*}
 is the Fourier-Mukai  kernel of the identity functor on $\Db(\op{Qcoh}^{\gl}Z)$.
\end{lem}
\begin{proof}
For an $R$-module $M$, the integral transform associated to $\mdx$ is given by 
\begin{equation*}
\Phi_{\mdx}(\widetilde{M}):=\left[{\bf R}\pi_{2*}\left(\mdx\otimes^{\bf{L}}{\bf L}\pi_1^*\widetilde{M}\right)\right]^{\gl},
\end{equation*}
where $\pi_i$ are the natural $\gl$-equivariant projections $Z\ktimes Z\ra Z$ (see \cite[Section 2]{BFKe} for background). Our desired result is a consequence of the following calculation:
\begin{align*}
\Phi_{\mdx}(\widetilde{M})&=\left[{\bf R}\pi_{2*}\left[\left((\pi_Z\ktimes\sigma_Z)_*\sch_{\Delta_Z}\right)\otimes^{\bf{L}}_{\sch_{Z\ktimes Z}}({\bf L}\pi_1^*\widetilde{M})\right]\right]^{\gl}\\
&\cong\left[\pi_{2^*}(\pi_Z\ktimes\sigma_Z)_*\left[\sch_{\Delta_Z}\otimes_{\sch_{\Delta_Z}}\left((\pi\ktimes\sigma_Z)^*\pi^*_1\widetilde{M}\right)\right]\right]^{\gl}\\
&\cong\left[ \sigma_{Z*}\pi_Z^*\widetilde{M}\right]^{\gl}\\
&\cong\left[\sch_{\Delta_Z}\otimes_{\sch_Z}\widetilde{M}\right]^{\gl}\\
&\cong\left[\left(\sch_{\gl}\otimes_\kk\sch_{Z}\right)\otimes_{\sch_Z}\widetilde{M}\right]^{\gl}\\
&\cong\left[\sch_{\gl}\otimes_\kk\widetilde{M}\right]^{\gl}\\
&\cong\widetilde{M},
\end{align*}
where the first isomorphism follows from the projection formula, and the last follows from Lemma \ref{lem:invar}. Furthermore, on the second isomorphism we may forego the process of deriving these functors as they are either exact or remain an adapted class (as discussed above).  As each of the above isomorphisms is induced by a natural transformation of functors, the above sequence yields a natural isomorphism between $\Phi_{\mdx}$ and the identity functor.
\end{proof}

\noindent
We now define the natural generalization of the object  $ Q $ from \cite[Defn 2.1.6]{BDF}.

\begin{defn}\label{defn:functor}
Given the scheme $Z=\spec(R)$ as before, define
 \begin{equation*}
Q_Z:=\left\lan \hpi_Y(R),\hsi_Y(R),C \right\ran\subseteq \kk\left[ \gl\ktimes Z\right].
\end{equation*}
that is the $\kk$-subalgebra of $\kk[\gl\times Z]$ generated by the images of $\sigma^\sharp_Z, \pi^\sharp_Z$ and the image of the inclusion $\kk[\op{End}(V)]\hookrightarrow\kk[\gl\ktimes Z]$. For ease of notation we denote $\sQ_Z:=\spec(Q_Z)$.
\end{defn}

\begin{rem}
Similar to the functor $Q$ in \cite[Def 2.1.6]{BDF} our definition provides a partial compactification of the action of $\gl$ on $Z$. For ease of reference we recall the definition of a partial compactification next.
\end{rem}

\begin{defn}
Let $G$ be an algebraic group and $Z$ a $\kk$-scheme with $G$-action. Also, let $\tilde{Z}$ be a $\kk$-scheme together an action of $G\times_\kk G$ which is equipped with a $(G\ktimes G)$-equivariant open immersion
\begin{equation*}
i:G\ktimes Z\hookrightarrow \widetilde{Z},
\end{equation*}
as well as a $(G\times_\kk G)$-equivariant morphism 
\begin{equation*}
(p,s):\widetilde{Z}\ra Z\times_\kk Z
\end{equation*}
such that the following diagram commutes
\begin{center}\begin{tikzcd}
&&\widetilde{Z}\ar[dd,shift right=.5ex,"p"']\ar[dd,shift left=1ex,"s"]\\
&\\
G\ktimes Z\ar[uurr,"i"] \arrow[rr, shift left=1ex,"\pi"] \arrow[rr,shift right=.5ex,"\sigma"']
&&Z
\end{tikzcd}\end{center}
where $\sigma$ is the action of $G$ on $Z$ and $\pi$ is the projection to $Z$. In this case, we refer to $\widetilde{Z}$, with the maps $p,s,i$, as a \emph{partial compactification of the action of} $G$ \emph{on} $Z$.
\end{defn}

\begin{ex}
In the case that $\dim V$ = 1, $ \gl = \mathbb{G}_m $, and  the definition of $Q$ given here recovers that found in \cite{BDF}.
\end{ex}

\begin{lem}
There are morphisms
\begin{center}\begin{tikzcd}
\sQ_Z\arrow[r, shift left=.5ex,"p"] \arrow[r, shift right=.5ex,"s"']&Z.
\end{tikzcd}\end{center}
which compose with the open immersion $ \Delta_Z \to \sQ_Z $ to give the morphisms  $\pi_Z$ and $\sigma_Z$.
\end{lem}
\begin{proof}
By definition, the maps $\hpi_Z$ and $\hsi_Z$ both have images which lie in $Q_Z$.
\end{proof}

\begin{lem}\label{lem:yeesh}
We have an isomorphism
 \begin{equation}\label{qiso}
Q_Z\cong\kk\left[A^L,B^L,A^R,B^R,C\right]/\left\lan B^L-B^RC,A^R-CA^L \right\ran \cong\kk[A^L,B^R,C].
\end{equation}
\end{lem}

\begin{proof}
We provide the reader with an easily verifiable isomorphism defined on the generators by
\begin{align*}
c_{ij}&\mapsto c_{ij},
\;\;\;\;\;\hpi_Z(a_{ij})\mapsto a_{ij}^L,
\;\;\;\;\;\hpi_Z(b_{ij})\mapsto b_{ij}^L,\\
\hsi_Z(a_{ij})&\mapsto a_{ij}^R,
\;\;\;\;\;\hsi_Z(b_{ij})\mapsto b_{ij}^R. \qedhere
\end{align*}

\end{proof}

\begin{rem}\label{itsye}
It follows from Equation~\eqref{qiso} that $\sQ_Z$ is  isomorphic to the closed subvariety of $(Z\ktimes Z)\ktimes\op{End}(V)$, consisting of the following points
\begin{equation*}
\left\{\left.(\psi_1,\psi_2,\psi_3,\psi_4,\varphi)\right| \psi_1=\psi_3\circ\varphi,\; \psi_4=\varphi\circ\psi_2 \right\}
\end{equation*}
\end{rem}

\begin{lem}
The scheme $\sQ_Z$ admits a $(\gl\ktimes\gl)$-action, denoted $\sigma_{\sQ_Z}$, which is uniquely defined  by the co-action 
\begin{equation*}
\hsi_{\sQ_Z}:\kk[A^L, B^R, C]\ra\kk[D^L , (\det D^L)^{-1}]\otimes\kk[D^R , (\det D^R)^{-1}] \otimes \kk[A^L, B^R, C],
\end{equation*}
which maps the generators 
 \begin{align*}
  b^R_{ij}&\mapsto \left(\det(D^R)\right)\inv\sum_{r=1}^n \op{Adj}(D^R)_{rj}\otimes_\kk b^R_{ir},\\
  a^L_{ij}&\mapsto \sum_{r=1}^n d^L_{ir}\otimes_\kk a^L_{rj},\\
  c_{ij}&\mapsto \left(\det(D^L)\right)\inv \sum_{r=1}^n \op{Adj}(D^L)_{sj} \otimes_\kk d^R_{ir}\otimes_\kk c_{rs},
\end{align*} where $ \op{Adj}(D) $ is the adjugate of the matrix $ D $.
\end{lem}

\begin{proof}
To check commutativity of the appropriate diagrams, we can pass to $\kk[G \times Z]$. Since $(Q_Z)_{\op{det} C} = \kk[G \times Z]$ and 
$\op{det} C$ is a non-zero-divisor, base change is flat. For $\kk[G \times Z]$, we are describing, in coordinates, the co-action 
corresponding to the previously specified $G \times G$ action on $\Delta_Z$. 
\end{proof}

The next lemma gives explicit descriptions of the two module structures that $Q_Z$ possesses. 

\begin{lem}\label{repQ}
For $Z=\spec(\kk[A,B])$ , we have the following two $\kk[A,B]$-module structures on $Q_Z$  given by $p^\sharp$ and $s^\sharp$, respectively:
\begin{align*}
p^\sharp : \kk[A,B] & \to \kk[A^L, B^R, C] \\
B&\mapsto B^RC \\
A&\mapsto A^L \\
s^\sharp : \kk[A,B] & \to \kk[A^L, B^R, C] \\
B&\mapsto B^R \\
A&\mapsto CA^L 
\end{align*}
\end{lem}
\begin{proof}
These are just the maps induced by the description of $Q_Z$ from Lemma \ref{lem:yeesh} under the identification
\[
Q_Z = \kk[A^L, B^R, C]. \qedhere
\]
\end{proof}

\begin{rem}\label{sorry}
We could define $ \sQ $ and $ \Delta $ as functors from affine varieties with a $ \gl $-action to affine schemes over $ \kk[\op{End}(V)] $ with a $ \gl \times \gl  $-action, thereby generalizing the work of \cite{BDF}.  However, as our focus in this paper is  the case of Grassmann flops as considered by \cite{DS}, we do not consider these generalizations.
\end{rem}

Now, we prove some properties of $ Q $ that will be used in Section~\ref{sec:bl} to prove the fullness of a Fourier-Mukai transform constructed using $ Q $.

\begin{lem}\label{notor}
For an object $Z=\spec(R)$ as before, we have
\begin{equation*}
\op{Tor}_i(_pQ_Z,\,_sQ_Z)=0
\end{equation*}
for all $i>0$, where the subscripts preceding $Q_Z$ denote the $R$-module structures given by $p^\sharp$ or $s^\sharp$, respectively. 
\end{lem}

\begin{proof}

Let $R:=\kk[A,B]$ as in Equation~\eqref{R}. By Lemma \ref{repQ}, we have
\begin{align*}
_pQ_Z&\cong \kk[A,B,B',C]\left/\left(B-B'C\right)\right. \\
_sQ_Z&\cong \kk[A,B,A',C]\left/\left(A-CA'\right)\right.
\end{align*}
Let us compute $ Q_Z \;_s\hspace{-1mm}\otimes_p^{\bf{L}} Q_Z $ using the above expressions:
\begin{align*}
Q_Z \;_s\hspace{-1mm}\otimes_p^{\bf{L}} Q_Z &= \kk[A,B,A',C]\left/\left(A-CA'\right)\right. \otimes^{\bf{L}}_{\kk[A,B]} \kk[A,B,B',C]\left/\left(B-B'C\right)\right. \\
&\cong \kk[A,B,A',C]\left/\left(A-CA'\right)\right. \otimes_{\kk[A,B]} \mathcal{K}_{\kk[A,B,B',C]}(B-B'C) \\
&\cong \mathcal{K}_{\kk[A,B,A',B',C_1,C_2]\left/\left(A-C_2A'\right)\right.}(B-B'C_1) \\
&\cong \mathcal{K}_{\kk[B,A',B',C_1,C_2]}(B-B'C_1),
\end{align*} where we resolved the regular sequence $\left(B-B'C\right)$ by the Koszul complex, denoted by $\mathcal{K}$, on the second line.

Finally, we see that the sequence $ (B-B'C_1) $ is still regular in the ring $ \kk[B,A',B',C_1,C_2] $ and hence all the higher homologies vanish.
\end{proof}

\begin{notation}\label{notation:Qtilde}
We denote by
\begin{equation*}
\widehat{Q}_Z:=(p\ktimes s)_*\mathcal O_{\sQ_Z}.
\end{equation*}
the sheaf of modules over $Z$ associated to $Q_Z$.
We will use the same notation in the derived setting (see  Section~\ref{section:derivedQ}, particularly Remark~\ref{rem:Qunderived}). Furthermore, as $\Delta_Z$ is an open subset of $\sQ_Z$ we will denote the natural open immersion as
\begin{equation*} \label{equation:eta}
\eta:\gl\ktimes Z\ra\sQ_Z.
\end{equation*}
\end{notation}

It will be useful to consider the following open covers of the quasi-affine sets $ U^+ $ and $ U^- $ defined before. Let
\begin{align}\label{uplus}
U^{+}&=\bigcup_{J\subseteq[m],\; |J|=d}U^{+}_J,
\end{align}
\begin{align}\label{uminus}
U^-=\bigcup_{I\subseteq[m'],\; |I|=d}U^{-}_I,
\end{align}
where
\begin{align*}
U^+_I&:=\spec\left(\kk\left[A,B,\left(\det(A_{[d],I}\right)\inv\right]\right)\\
U^-_J&:=\spec\left(\kk\left[A,B,\left(\det(B_{J,[d]}\right)\inv\right]\right)
\end{align*}
and (for example) $\det\left(A_{[d],I} \right)$ denotes the $(d\times d)$ minor of $A$ consisting of the rows indexed by $I$. Therefore, we have the following affine open covers:
\begin{equation}\label{equation:cover1}
U^+\times_{\rox} U^-=\bigcup_{I\subseteq[m],\;\; J\subseteq[m'],\;\;|I|=|J|=d}U^+_I\times_{\rox} U^-_J,
\end{equation}

\begin{equation}\label{equation:cover2}
U^+\times_{\kk} U^-=\bigcup_{I\subseteq[m],\;\; J\subset[m'],\;\;|I|=|J|=d}U^+_I\times_{\kk} U^-_J,
\end{equation}
where $\rox := \Spec(\kk[A, B]^{\text{GL}(V)})$ denotes the invariant theoretic quotient of $Z$.

\begin{lem}\label{lemma:cong}
There is an isomorphism
\begin{equation*}
\kk\left[Z\times_{\rox}Z\right]\cong \kk\left[A^L,B^L,A^R,B^R\right]/\left\lan B^LA^L-B^RA^R \right\ran , 
\end{equation*}
where the generators and relations are as in Definition \ref{defn:functor}.
\end{lem}

\begin{proof}
 From Weyl's fundamental theorems for the action of $\gl$ (for example see \cite[Chapter 2.1]{kraft} or the original text \cite{wey46}) we have
 \begin{displaymath}
  \rox = \{ D \in \operatorname{Hom}_{\kk}(W,W) \mid \operatorname{rank} D \leq \operatorname{dim} V \}.
 \end{displaymath}
The map $Z \to \rox$ is thus given by the homomorphism 
 \begin{align*}
  \kk[\rox] & \to k[A,B] \\
  D & \mapsto BA. 
 \end{align*}
 Hence, 
 \begin{displaymath}
  \kk\left[Z\times_{\rox}Z\right] = \kk[A^L,B^L] \otimes_{\kk[\rox]} \kk[A^R,B^R] \cong \kk[A^L,B^L,A^R,B^R]/(B^LA^L - B^RA^R). \qedhere
 \end{displaymath}
\end{proof}

\begin{lem}
There exists a morphism
\begin{displaymath}
 \kappa := p^\# \otimes s^\# : \kk[A^L,B^L] \otimes_{\kk[\rox]} \kk[A^R,B^R] \to Q_Z. 
\end{displaymath}
\end{lem}
\begin{proof}
This follows since $p^\#$ and $s^\#$ are equal on $\kk[\rox]$, by definition.
\end{proof}

\begin{lem}\label{lemma:contains}
With the conventions above we have the following containment of ideals in the ring $\kk[A^L,B^L,A^R,B^R,C]$:
 \begin{equation*}
 \left\lan B^LA^L-B^RA^R\right\ran\subset\left\lan B^L-B^RC,A^R-CA^L \right\ran.
\end{equation*}
 \end{lem}
 
\begin{proof}
 This follows from 
 \begin{displaymath}
  (B^L-B^RC)A^L + B^R(CA^L - A^R) = B^LA^L - B^RA^R. \qedhere
 \end{displaymath}
\end{proof}

\begin{prop}\label{localiso}
Let $Z=\spec(\kk[A,B])$ as usual and $Q_Z$ as in Equation~\eqref{qiso}. Let $\widehat{Q}_Z|_{U^+\ktimes U^-}$ be the restriction of $\widehat{Q}_Z$ to the open subset $U^+\ktimes U^-\subset Z\ktimes Z$. Then $\kappa$ restricts to an isomorphism 
\begin{equation*}
 \kappa|_{U^+ \ktimes U^-} : \left.\widehat{Q}_Z\right|_{U^+\ktimes U^-}\overset{\sim}{\longrightarrow} \sch_{U^+\times_{\rox}U^-}.
\end{equation*}
\end{prop}

\begin{proof}
We look affine-locally using the covers of Equations \ref{equation:cover1} and \ref{equation:cover2}. We need only show that under the above localization the map $\kappa:\kk[Z\times_{\rox}Z]\ra Q_Z$ becomes an isomorphism. For surjectivity, it suffices to show that there is an element (we find two such) which map to $C$. Indeed, we have
\begin{align*}
\Big((B^R)_{J[d]}\Big)\inv B^L&\mapsto C\\
 A^R\Big((A^L)_{[d]I}\Big)\inv&\mapsto C,
\end{align*}
easily verified by the relations $B^L-B^RC$ and $A^R-CA^L$ in $Q(\kk[A,B])$ given in Definition \ref{defn:functor}. 

For injectivity, it suffices to check that under this localization we have the containment
\begin{equation*}
\left\lan B^L-B^RC,A^R-CA^L \right\ran\subset \left\lan B^LA^L-B^RA^R\right\ran, 
\end{equation*}
since the opposite containment is Lemma~\ref{localiso}. To see this, simply note that by multiplying by the appropriate elements in the above identification, we have
\begin{equation*}\label{1}
 (B^L)_{J[d]}=(B^R)_{J[d]} C\;\;\;\;\;\text{ and }\;\;\;\;\;(A^R)_{[d]I}= C(A^L)_{[d]I}.
\end{equation*}
Hence, multiplying by the appropriate units in our localization, we have
\begin{equation*}
\left\lan B^L-B^RC,A^R-CA^L\right\ran=\left\lan (B^R)_{J[d]}\left(A^R-CA^L\right), \left( B^L-B^RC\right)(A^L)_{[d]I} \right\ran .
\end{equation*}
For example, by Equation~\eqref{1}, we have
\begin{align*}
 (B^R)_{J[d]}\left(A^R-CA^L\right)= \left((B^R)_{J[d]}A^R-(B^L)_{J[d]}A^L\right)\in\left\lan B^LA^L-B^RA^R\right\ran,
\end{align*}
while the other relation follows similarly. This gives our desired isomorphism.
\end{proof}

Consider the restriction $Q|_{U^+ \times U^-}$. By descent, we have a corresponding object $P$ on the quotient $Z^+ \times Z^-$.

\begin{thm} \label{theorem:q_fib_prod}
 We have an isomorphism 
 \begin{displaymath}
  P \cong \mathcal O_{Z^+ \times_{Z_0} Z^-}
 \end{displaymath}
\end{thm}

\begin{proof}
 This follows immediately by passing to the quotient in Proposition~\ref{localiso}. 
\end{proof}

We now examine a useful invariant when studying kernels in the next subsection. Note that for $Z$,  the tensor product $\qsquare$ is equipped with a natural $\gl^{\times 4}$-action. This induces a $\gl^{\times 3}$-action, which we denote
\begin{equation*}
\sigma_3:\gl ^{\times 3} \ktimes\sQ_Z^{\times 2}\ra\sQ_Z^{\times 2}
\end{equation*}
and is defined as the product of the following compositions
\begin{center}\begin{tikzcd}
&\gl^{\times 3}\ktimes\sQ_Z^{\times 2}\ar[dl,"\pi_{1,2,4}"']\ar[dr,"\pi_{2,3,5}"]&\\
\gl^{\times 2}\ktimes\sQ_Z\ar[d,"\sigma_{\sQ_Z}"']&&\gl^{\times 2}\ktimes\sQ_Z\ar[d,"\sigma_{\sQ_Z}"]\\
\sQ_Z&&\sQ_Z
\end{tikzcd}\end{center}
Here $\pi_{i,j,k}:\gl^{\times 3} \ktimes\sQ_Z^{\times 2} \to \gl ^{\times 2} \ktimes \sQ_Z^{}$ is the projection onto the $i^{\text{th}}$, $j^{\text{th}}$ and $k^{\text{th}}$ components. For any ring $T$ with $\gl ^{\times 3}$-action, we will denote the invariant subring  associated to the action corresponding to the middle component of $\gl^{\times 3}$ by $T^{\bowtie}$. The notation $(-)^{\bowtie}$ is suggestive of pinching a module in the middle. Since taking invariants is functorial for equivariant morphisms, we obtain the following: 

\begin{lem}\label{lem:Qpinch}
The following diagram commutes
\begin{center}\begin{tikzcd}
&\left(Q_Z\hspace{1mm}_p\hspace{-1mm}\otimes_\sigma \Delta_Z\right)^{\bowtie}\ar[dr,"\sim"]&\\
\left(Q_Z\hspace{1mm}_p\hspace{-1mm}\otimes_s Q_Z\right)^{\bowtie}\ar[ur,"(1\otimes\eta)^{\bowtie}"]\ar[dr,"(\eta\otimes 1)^{\bowtie}"']&&Q_Z\\
&\left(\Delta_Z\hspace{1mm}_\pi\hspace{-1mm}\otimes_sQ_Z\right)^{\bowtie}\ar[ur,"\sim"']&
\end{tikzcd}\end{center}
Furthermore, the morphism $\rho_Z:\left(Q_Z\hspace{1mm}_p\hspace{-1mm}\otimes_sQ_Z\right)^{\bowtie}\ra Q_Z$ is an isomorphism. 
\end{lem}

\begin{proof}
First recall that we have a presentation from Lemma \ref{repQ} of $_sQ_Z$ and $_pQ_Z$, which for ease of calculation we set the following simplified notation, with the hope that no confusion arises:
\begin{align*}
_pQ_Z&\cong \kk[A,B^L,B^R,C]/\left(B^L-B^RC\right)\cong \kk[A,B^R,C]\\
&:=\kk[A,B,C]\\
_sQ_Z&\cong \kk[A^L,A^R,B,C]/\left(A^R-CA^L\right)\cong\kk[A^L,B,C]\\
&:=\kk[A,B,C]
\end{align*}
Further we recall the notational preference that for $\kk$-algebras $R,S$ and $r\in R$, $s\in S$ that the following pure tensors will be denoted: $r\otimes1:=r^L$ and $1\otimes s:=s^R$. With these conventions we have the following presentations of rings:
\begin{align*}
Q_Z\hspace{1mm}_p\hspace{-.3mm}\otimes_sQ_Z&\cong \kk[A^L,A^R,B^L,B^R,C^L,C^R]\left/\left(B^LC^L-B^R,A^L-C^RA^R\right)\right.\\
&\cong\kk[A^R,B^L,C^L,C^R]\\
\kk[\Delta_Z]\hspace{1mm}_\pi\hspace{-.3mm}\otimes_sQ_Z&\cong \kk[A^L,A^R,B^L,B^R,C^L,C^R,\det(C^L)\inv]\left/\left(B^L-B^R,A^L-C^RA^R\right)\right.\\
&\cong\kk[A^R,B^L,C^L,C^R,\det(C^L)\inv]\\
Q_Z\hspace{1mm}_p\hspace{-.3mm}\otimes_\sigma\kk[\Delta_Z]&\cong\kk[A^L,A^R,B^L,B^R,C^L,C^R,\det(C^R)\inv]
\left/\left(B^L-B^R\left({C^R}\right)
\inv,A^L-C^RA^R\right)\right.\\
&\cong\kk[A^R,B^L,C^L,C^R,\det(C^R)\inv]
\end{align*}
Hence, commutativity of the above diagram is clear. Furthermore, one verifies that we have an isomorphism $\kk[\Delta_Z]\hspace{1mm}_\pi\hspace{-.3mm}\otimes_sQ_Z\cong Q_Z\hspace{1mm}_p\hspace{-.3mm}\otimes_\sigma\kk[\Delta_Z]$, and thus 
\begin{equation*}
\left(\kk[\Delta_Z]\hspace{1mm}_\pi\hspace{-.3mm}\otimes_sQ_Z\right)^{\bowtie}\cong \left(Q_Z\hspace{1mm}_p\hspace{-.3mm}\otimes_\sigma\kk[\Delta_Z]\right)^{\bowtie}.
\end{equation*}
 It is clear that the maps on the right-hand side of the diagram are  isomorphisms since $\kk[\Delta_Z]$ is the kernel of the identity by Lemma \ref{kernelid}. We claim that
\begin{align*}
\left(Q_Z\hspace{1mm}_p\hspace{-.3mm}\otimes_sQ_Z\right)^{\bowtie}&=\kk[A^R,B^L,C^L\cdot C^R]\\
\left(Q_Z\hspace{1mm}_p\hspace{-.3mm}\otimes_\sigma \kk[\Delta_Z]\right)^{\bowtie}&=\kk[A^R,B^L,C^L\cdot C^R]
\end{align*}
from which it follows that these rings are isomorphic. This claim is simply Weyl's Theorem for the invariants of $\kk[V\otimes V^{\vee}]$.
\end{proof}



\subsection{The integral kernel}\label{section:derivedQ}
We now use $\sQ$ to construct Fourier-Mukai kernels. We begin by recalling the following from \cite[Definition 3.1.4]{BDF}.
\begin{defn}
Let $\widetilde{Z}$ be a partial compactification of an action $\sigma:G\ktimes Z\ra Z$, with maps $p, s,$ and $i$ as above. We define the \emph{boundary} of $\widetilde{Z}$ to be
\begin{equation*}
\partial^s_{\tilde{Z}}:=\widetilde{Z}\setminus i\left(G\ktimes Z\right),
\end{equation*}
the $s$-\emph{unstable locus} to be
\begin{equation*}
Z_s^{\op{us}}:=s\left(\partial_{\widetilde{Z}}\right),
\end{equation*}
and the $s$-\emph{semistable locus} to be
\begin{equation*}
Z_s^{\op{ss}}:=Z\setminus Z^{\op{us}}.
\end{equation*}
One similarly defines  the $p$-\emph{unstable} and $p$-\emph{semistable} loci.
\end{defn}
\begin{rem}
It follows from \cite[Example 3.1.10]{BDF} that  the $s$-semistable locus $Z_s^{\op{ss}}$ coincides with $U^+$ from Equation~\eqref{uplus}. Similarly, the $p$-semistable locus $Z_p^{\op{ss}}$ coincides with $U^-$ from Equation~\eqref{uminus}.
\end{rem}

\begin{defn}
For an object $Z$ as before, we let
\begin{equation*}
\pQ_Z:=(p\times s)_*\sch_{\sQ_Z}\in\Db\left(\op{Qcoh}^{\gl\ktimes\gl}Z\ktimes Z\right),
\end{equation*}
where the pushforward is understood to be derived. We denote by $\pQs_Z$ the quasi-coherent sheaf on $Z_s^{\op{ss}}\ktimes Z$ realized by restricting $\pQ_Z$ from $Z\ktimes Z$. That is,
\begin{equation*}
\pQs_Z=(j\times 1_Z)^*\pQ_Z,
\end{equation*}
where $j:Z_s^{\op{ss}}\ra Z$ is the inclusion. Finally, taking $\pQs_Z$ as the Fourier-Mukai kernel, we have the functor
\begin{equation}\label{eq:qplus}
\Phi_{\pQs_Z}:\Db\left(\op{Qcoh}^{\gl}Z_s^{\op{ss}}\right)\ra\Db\left(\op{Qcoh}^{\gl}Z\right).
\end{equation}
\end{defn}
\begin{rem}\label{rem:Qunderived}
Since the functor $(p\times s)_*$ is exact, $\pQ_Z$ is just the $\gl$-linearized sheaf associated to $Q_Z$ with its $(p,s)$-bimodule structure given in Lemma \ref{lemma:dubact}. This justifies our use of $\widehat{Q}_Z$ in Notation~\ref{notation:Qtilde}. 
\end{rem}

\begin{lem}
 The functor $\Phi_{\pQs_Z}$ is faithful.
\end{lem}
\begin{proof}
Our proof follows from the fact that the functor 
\begin{equation*}
i^*:\Db(\op{Qcoh}_{\gl}(Z^{ss}_s))\ra\Db(\op{Qcoh}_{\gl}(Z))
\end{equation*}
 is the left inverse of $\Phi_{\widehat{Q}^+_{Z}}$. To see this, note that for any maximal minor $m$ of $B$, we have $R_m\otimes_sQ_Z\cong\kk[\gl]\kten R=\kk[\Delta_Z]$. Indeed, inverting a minor on the left amounts to inverting the determinant of $C$. Since $\Delta_Z$ is the kernel of the identity, we obtain the desired result. 
\end{proof} 

The fullness of this functor depends on certain localization properties, which are the focus of the next section. 



\subsection{Bousfield localizations}\label{sec:bl}
This section recalls Bousfield (co)-localizations which will be used to establish fullness of the functor $\Phi_{\pQs}$ from Equation~\eqref{eq:qplus}. We recall that the existence of a Bousfield triangle produces a semi-orthogonal decomposition, and we show that the essential image of our functor is an inclusion into one of these pieces. We refer the reader to \cite{kra} for a more detailed treatment of these concepts. While the proofs of the statements refer to \cite{BDF} we recall all of the statements here for ease of reference.

\begin{defn}
Let $\mathcal{T}$ be a triangulated category. A \emph{Bousfield localization} is an exact endofunctor $L:\mathcal{T}\ra\mathcal{T}$ equipped with a natural transformation $\delta:1_{\mathcal{T}}\ra L$ such that:
\begin{enumerate}
\item[a)] $L\delta=\delta L$ and
\item[b)] $L\delta:L\ra L^2$ is invertible.
\end{enumerate}
A \emph{Bousfield co-localization} is given by an endofunctor $C:\mathcal{T}\ra\mathcal{T}$ equipped with a natural transformation $\epsilon:C\ra 1_{\mathcal{T}}$ such that:
\begin{enumerate}
\item[a)] $C\epsilon=\epsilon C$ and 
\item[b)] $C\epsilon:C^2\ra C$ is invertible.
\end{enumerate}
\end{defn}
\begin{defn}\label{triangle}
Assume there are natural transformations of endofunctors
\begin{equation*}
C\overset{\epsilon}\ra1_{\mathcal{T}}\overset{\delta}\ra L
\end{equation*}
of a triangulated category $\mathcal{T}$ such that
\begin{equation*}
Cx\overset{\epsilon_{Cx}}\longrightarrow x\overset{\delta_x}\longrightarrow Lx
\end{equation*}
is an exact triangle for any object $x$ of $\mathcal{T}$. Then we refer to $C\ra1_{\mathcal{T}}\ra L$ as a \emph{Bousfield triangle} for $\mathcal{T}$ when any of the following equivalent conditions are satisfied:
\begin{enumerate}
\item[1)] $L$ is a Bousfield localization and $C(\epsilon_x)=\epsilon_{C_x}$
\item[2)] $C$ is a Bousfield co-localization and $L(\delta_x)=\delta_{L_x}$
\item[3)] $L$ is a Bousfield localization and $C$ is a Bousfield co-localization.
\end{enumerate}
\end{defn}

\noindent
For a proof that the above properties are indeed equivalent, we refer the reader to \cite[Definition 3.33]{BDF}. Denoting $S:=\kk[\Delta_Z]/Q_Z$, we have morphisms 
\begin{align*}
Q_Z\overset{\eta^\sharp}\ra\kk[\Delta_Z]\ra S\ra Q_Z[1]
\end{align*}
in $\Db(\op{Mod}^{\gl} \kk[Z])$, where $\eta^\sharp$ is the morphism induced by $\eta$ as in Equation~\ref{equation:eta}. This yields an exact triangle. Furthermore, if we let $\hat{\eta}:\Phi_{\pQ_Z}\ra1$ denote the morphism induced by $\eta^\sharp$, we see that for any $x$ in $\Db(\op{Qcoh}^{\gl}Z)$ the following is also exact:
\begin{equation*}
\Phi_{\pQ_Z}(x)\ra x\ra\Phi_{\widetilde{S}}(x)
\end{equation*}
With these observations in mind, we present one of the main results of this section.
\begin{prop}\label{prop:triangle}
The triangle of functors 
\begin{equation*}
\Phi_{\pQ_Z}\overset{\hat{\eta}}\ra1\ra\Phi_{\widetilde{S}}.
\end{equation*}
is a Bousfield triangle.
\end{prop}

\begin{proof}
This follows identically as in \cite[Lemma 3.3.6]{BDF}, by Lemma \ref{notor} and Lemma \ref{lem:Qpinch}.  
\end{proof}

We are now ready to prove that $\Phi_{\pQs}$ is full. Let $J_+:=j_*\circ j^*$, where $j:Z_s^{\op{ss}}\ra Z$ is the natural inclusion, and let $\Gamma_+$ be the local cohomology.

\begin{prop}\label{phidecomp}
There is a semi-orthogonal decomposition 
\begin{equation*}
\mathsf{D}(\op{Qcoh}^{\gl}Z)=\langle\op{Im}\Phi_{\widetilde{S}},\op{Im}\Phi_{\pQs},\op{Im}\Phi_Q\circ\Gamma_+\rangle,
\end{equation*}
where $\op{Im}$ denotes the essential image. Furthermore, $\Phi_{\pQs}$ is fully-faithful.
\end{prop}

\begin{proof}
This follows identically to the proof of Proposition 3.3.9 in \cite{BDF}
\end{proof}

Letting $j':Z^{ss}_p\ra Z$ be the inclusion and $\Gamma_+$ its local cohomology, we have the following dual statement.
\begin{prop}\label{phidecompminus}
There is a semi-orthogonal decomposition 
\begin{equation*}
\mathsf{D}(\op{Qcoh}^{\gl}Z)=\langle\op{Im}\Phi_{\widetilde{S}},\op{Im}\Phi_{\widehat{Q}^-},\op{Im}\Phi_Q\circ\Gamma_-\rangle,
\end{equation*}
where $\op{Im}$ denotes the essential image. Furthermore, $\Phi_{\widehat{Q}^-}$ is fully-faithful.
\end{prop}



\section{A geometric resolution}\label{geotech}
For this section, we will denote $Z$ as the scheme
\begin{equation*}
\hom(V,W) \oplus \hom(W',V).
\end{equation*}
Having established that $\Phi_{\pQs}$ is fully faithful, the remaining objective of this work is to examine the essential image of the functor $\Phi_{\pQs}$. We will show that this image is generated by an exceptional collection first discovered by Kapranov in \cite{Kapp}.  The method which we use is based on the underlying techniques of the well known `geometric technique' of Kempf (see e.g. \cite{wey}). 

\subsection{A sketch of Kempf}
The objective of the method of Kempf is to provide a free resolution of special modules by pulling back to a trivial geometric bundle over a projective variety. 

Consider an algebraic variety $Y$. The total space of the sheaf $ \mathcal{O}_Y^{\oplus n} $ is the scheme $ Y \times \mathbb{A}^n$. Now let $X$ be the total space of a locally free sheaf  $\mathcal{F} \subset \mathcal{O}_Y^{\oplus n}$ on $ Y $.  Let $ \pi $ denote the projection $ Y \times \mathbb{A}^n\ra Y$.

We have the exact sequence of locally free sheaves on $ Y \times \mathbb{A}^n$
\begin{center}\begin{tikzcd}
0\ar[r]& \pi^*\mathcal{F}\ar[r]& \pi^*\mathcal{O}_{ Y}^{\oplus n}\ar[r,"f"]&\pi^*\mathcal{T}\ar[r]&0,
\end{tikzcd}\end{center} where $ \mathcal{T} $ is the quotient sheaf.

Consider the section $ s:= f \circ \op{taut} : \mathcal{O}_{Y \times \mathbb{A}^n} \to \pi^* \mathcal{T} $, where $ \op{taut} $ denotes the tautological section of $ \pi^*\mathcal{O}_{ Y}^{\oplus n} $ on $ Y \times \mathbb{A}^n $. Then, we have the following statement.

\begin{prop}\label{kempf}
	With the above notation, a locally free resolution of the sheaf $\sch_{X}$ as a $\sch_{ Y \times \mathbb{A}^n}$-module is given by the Koszul complex
	\begin{equation*}
	\mathcal{K}\left(s\right)_{\bullet}: 0\ra\bigwedge\nolimits^{\op{rnk}(\mathcal{T})}\left(\pi^*\mathcal{T}^\vee\right)\ra\hdots\ra\bigwedge\nolimits^2\left(\pi^*\mathcal{T}^\vee\right)\ra \pi^*\mathcal{T}^\vee\ra\sch_{ Y \times \mathbb{A}^n}
	\end{equation*}
	\end{prop}

\begin{proof}
	On the vanishing locus $ Z(s) $, the tautological section $ \op{taut} $ factors through $ \pi^*\mathcal{F} $. Hence, the vanishing locus is the total space of the sheaf $ \mathcal{F} $, which is $ X $. We see that the section is regular as the codimension of $ Z(s) $ equals the rank of the sheaf $ \pi^* \mathcal{T} $; and the Koszul complex resolves $ \mathcal{O}_X $. For more details, see \cite[Proposition 3.3.2]{wey}.
\end{proof}

\subsection{The resolution}\label{theresolver}
Now we are ready to present a resolution which will open a window to view $\op{Im}(\Phi_{\pQs})$. First recall that we set $\dim(V):=d$.  We define $\mathbf{Q}^+_Z$ as the base change:
\begin{center}\begin{tikzcd}
		\mathbf{Q}^+_Z \ar[r,]\ar[d]&\mathbf{Q}_Z\ar[d,"p \times s"]  \\
		Z^{\op{ss}}_s\times Z \ar[r,hook]&Z \times Z
\end{tikzcd}\end{center}

Let $\mathcal S$ be the tautological bundle on $\op{Gr}(d,W)$ i.e. the locally free sheaf on $\op{Gr}(d,W) = [\op{Hom}(V,W)^{\op{ss}}_s / \gl^L]$ corresponding to the $\gl^L$-representation $V$.  Then, we have the Euler sequence for the Grassmannian $ \op{Gr}(d,W) $:
\[
0 \to \mathcal S \to W \to \mathcal Q \to 0.
\]
Consider the pullback of the above sequence to $ \op{Gr}(d,W) \times \op{Hom}(W',V) $ along $ q $ and  apply $ \sheafhom(t^*V,-) $, where $ q : \op{Gr}(d,W) \times \op{Hom}(W',V) \to \op{Gr}(d,W) $  and $ t : \op{Gr}(d,W) \times \op{Hom}(W',V) \to \op{Hom}(W',V) $ are projections
\begin{equation}
0 \to \sheafhom(t^*V, q^*\mathcal S) \xrightarrow{\rho} \sheafhom(t^*V, q^*W) \xrightarrow{\Xi} \sheafhom(t^*V, q^*\mathcal Q) \to 0.
\label{eq: ses}
\end{equation}

Let us denote $ \mathcal{T} := \sheafhom(t^*V, q^*\mathcal{Q}) $. We denote the total space of the locally free sheaf $ \sheafhom(A, B)  $ as $ \mathbf{Hom}(A,B) $. From the discussion in the previous subsection, we get the following result:

\begin{lem}\label{koszul01}
 The following Koszul complex is a free resolution for $\sch_{\mathbf{Hom}(t^*V,q^*\mathcal{S})}$ as an $\sch_{\op{Gr}(d,W) \times Z}$-module.
\begin{equation}\label{koszul1}
\mathcal{K}\left(s\right)_{\bullet}: \bigwedge^{d(m-d)}\pi^* \mathcal{T}^\vee\ra \hdots \ra\bigwedge^2 \pi^*\mathcal{T}^\vee \ra \pi^*\mathcal{T}^\vee\ra\sch_{ \op{Gr}(d,W) \times Z}
\end{equation} where $ \pi : \op{Gr}(d,W) \times \Hom(V,W) \times \Hom(W',V) \to \op{Gr}(d,W) \times \Hom(W',V)  $ is the projection morphism.
\end{lem}
\begin{proof}
	We choose $ Y =  \op{Gr}(d,W) \times \Hom(W',V) $, and $ \mathcal{F} = \sheafhom(t^*V, q^*\mathcal S)$, and apply Proposition \ref{kempf}. Notice that the total space of $ \sheafhom(t^*V, q^*W) $ on $ \op{Gr}(d,W) \times \op{Hom}(W',V)  $ is $ \op{Gr}(d,W) \times Z $.
\end{proof}

Now, we can identify $[\mathbf{Q}_Z^+/\gl^L]$ as the total space $ \mathbf{Hom}(t^*V, q^*\mathcal S) $.

\begin{lem} \label{Q is a bundle}
	The quotient space $[\mathbf{Q}_Z^+/ \gl^L]$ is $\gl^R$-equivariantly isomorphic to  the total space $\mathbf{Hom}(t^*V, q^*\mathcal S)$ as schemes over $\op{Gr}(d,W) \times \op{Hom}(W',V)$.
\end{lem}
\begin{proof}
	
	Recall from Equation~\eqref{qiso}, that $\mathbf{Q}_Z$ is associated to the module
	\[
	\kk[A^L,B^R,C]
	\]
	Geometrically, we may view $\mathbf{Q}_Z$ as the total space of the locally free sheaf $\op{End}(V)$ over $\op{Spec}\kk[A^L, B^R]$.	Once we base change to the semistable locus and take the quotient with respect to the  $ \gl^L $ action, we get that $[\mathbf{Q}_Z^+/\operatorname{GL}(V)^L]$ is isomorphic to
	the total space 
	\begin{displaymath}
	\mathbf{Hom}(t^*V, q^*\mathcal S) \to \op{Gr}(d,W) \times \op{Hom}(W',V).
	\end{displaymath}

	Moreover, the inclusion, $\sheafhom(t^*V, q^*\mathcal S) \to \sheafhom(t^*V, q^*W)$ realizes it as a subspace of the total space $\mathbf{Hom}(t^*V, q^*W)$ over $ \op{Gr}(d,W) \times \op{Hom}(W',V)$ which is  $Z \times  \op{Gr}(d,W)$.  
	
	This inclusion $ \sheafhom(t^*V, q^*\mathcal S) \to \sheafhom(t^*V, q^* W) $ is induced by the ring homomorphisn
	\begin{align*}
	k[A^L,A^R,B^R] & \to k[A^L,B^R,C] \\
	A^L & \mapsto A^L \\
	A^R & \mapsto CA^L \\
	B^R & \mapsto B^R. 
	\end{align*}
	which is equivariant with respect to the remaining $\gl^R$-action.
\end{proof}

We denote $\pi_1 : [Z^{\operatorname{ss}}_s/\gl^L] \to [\operatorname{Hom}(V,W)^{\operatorname{ss}}_s/\gl^L]$ as the projection. Putting Lemma \ref{koszul01} and Lemma \ref{Q is a bundle} together, we get a resolution of the sheaf $( \pi_1 \times \op{Id}_Z)_* \pQs_{\bZ}$.
\begin{cor} \label{resolution}
The Koszul complex (\ref{koszul1}) is a locally free resolution of the sheaf $(\pi_1 \times \op{Id}_Z)_* \pQs_{\bZ}$ of  $\O_{ \op{Gr}(d,W) \times Z}$-modules.
\end{cor}

\begin{rem}
	We note that we could also have constructed a locally free resolution of $ \pQs_{\bZ} $ on $ Z^{\op{ss}}_s \times Z $ by the same method, and this will also lead to a similar proof as in the remainder of this paper.
\end{rem}



\section{Analyzing the integral transform}
\noindent
In this section, we show that the kernel $ \pQs_Z $ induces a derived equivalence for a Grassmann flop. We begin by showing that the essential image of this functor coincides with the  `window' description studied by Donovan and Segal in \cite[Section 3.1]{DS}.  Specifically, we will show that the image of $\Phi_{\pQs_{\bZ}}$ is generated by a collection of vector bundles corresponding to representations identified by Kapranov \cite{Kapp}.  

Let us recall Kapranov's collection. Consider the standard $ \gl $ representation $ V $, where $ GL(V) $ acts by left multiplication. Consider the Schur modules of $ V $ associated to a Young diagram (or equivalently, partition) $ \alpha $, and denote them by $ L_\alpha V $. Kapranov's collection  is defined by
\begin{equation*}
\kapreezer_{d,m}:=\Big\{L_{\alpha}V\;\Big|\;\alpha\in \text{Young diagrams of height $ \leq m-d $ and width $ \leq d $}\Big\}.
\end{equation*}

We also consider pull backs of these representations to $ \op{Gr}(d,W) $ along the structure morphism. As $ V $ pulls back to the tautological bundle $ \mathcal{S} $, the Schur functors $ L_{\alpha}V $ pull back to $  L_{\alpha}\mathcal{S} $ and these are the locally free sheaves considered by Kapranov. By abuse of notation, we will consider $ \kapreezer_{d,m} $ as a collection of locally free sheaves on $\hom(V,W)\oplus\hom(W',V)$ or $\hom(W',V)$ (again, by pulling back along the structure morphism). Note that when $\kk=\mathbb{C}$, this is exactly the dual of the zero$^{\text{th}}$ window $W_0$ from \cite[Section 3.1]{DS}. 

It is the objective of this section to show that the thick triangulated subcategory generated by elements of $\kapreezer_{d,m}$ is equivalent to $\im\left(\Phi_{\pQs_{\bZ}}\right)$. We show one containment in Proposition~\ref{container1}, which relies on the work of Section \ref{geotech}.

\subsection{Windows from a resolution}
Consider the projection $\pi_1 : Z^{\operatorname{ss}}_s \to \operatorname{Hom}(V,W)^{\operatorname{ss}}_s$.  
To demonstrate that the image of $\Phi_{\pQs_{\bZ}}$ is contained in $
\langle \kapreezer_{d,m} \rangle$, we exhibit a particular $\operatorname{GL}(V)^L \times \operatorname{GL}(V)^R$-equivariant resolution $\mathcal K_\bullet$  of $(\op{Id}_Z \times \pi_1) _* \pQs_{\bZ}$ over $\operatorname{Hom}(V,W)^{\operatorname{ss}}_s \times Z$ i.e.\ we resolve the kernel of the functor $\Phi_{\pQs_{\bZ}} \circ \pi_1^*$.  Equivalently, this is a $\operatorname{GL}(V)^R$-equivariant resolution of  $(\op{Id}_Z  \times \pi_1)_* \pQs_{\bZ}$ over $\op{Gr}(d,W) \times Z$. The resolution obtained in equation~\eqref{koszul01} in Section \ref{theresolver} is the one we are looking for. 

 In this subsection, we will show that the components $\mathcal K^i$ of the resolution have a filtration whose associated graded pieces are of the form $J \boxtimes K$ with $K \in \kapreezer_{d,m}$.  This decomposition of the Fourier-Mukai transform $\Phi_{\pQs_{\bZ}} \circ \pi_1^*$ yields a functorial way to describe $\Phi_{\pQs_{\bZ}} \circ \pi_1^*(M)$ using objects of $\kapreezer_{d,m}$ for all objects $\pi_1^*(M) \in \Db( [Z^{\op{ss}}_s / \operatorname{GL}(V)] )$.  As such objects generate $\Db( [Z^{\op{ss}}_s / \operatorname{GL}(V)] )$ this is enough to conclude the goal of this section,  $\im\left(\Phi_{\pQs_{\bZ}}\right) \subseteq\langle\lweezer_{d,m}\rangle$.

\begin{prop}\label{container1}
With notation as above, we have 
\begin{equation*}
\im\left(\Phi_{\pQs_{\bZ}}\right)\subseteq\langle\lweezer_{d,m}\rangle , 
\end{equation*}
where $\langle\lweezer_{d,m}\rangle$ is the thick triangulated subcategory generated by elements in $\lweezer_{d,m}$.
\end{prop}

\begin{proof}
By Corollary~\ref{resolution}, we have a quasi-isomorphism with the Koszul complex
\[
\mathcal{K}_\bullet \cong 
(\op{Id}_Z  \times \pi_1)_* \pQs_{\bZ}
\]
The components of the Koszul complex are $\bigwedge^l \pi^\ast \sheafhom(t^*V, q^*\mathcal Q)^\vee$ for $0 \leq l \leq d$.  We can appeal to the Cauchy Formula, e.g. \cite[Theorem 2.3.2(a)]{wey}, to get a filtration on $\bigwedge\nolimits^{\!i} t^\ast \pi^\ast \sheafhom(t^*V, q^*\mathcal Q)^\vee$ whose associated graded pieces are 
\begin{displaymath}
 \pi^\ast \left(\bigoplus_{|\lambda| = i} L_{\lambda} V \boxtimes L_{\lambda'} \mathcal Q^\vee \right).
\end{displaymath}
Thus, each term in the Koszul complex can be generated using iterated exact sequences from the locally free sheaves
\begin{displaymath}
 \pi^\ast \left( L_{\lambda} V \boxtimes L_{\lambda'} \mathcal Q^\vee \right). 
\end{displaymath}
These components, in turn, generate $\pQs_{\bZ}$.  Hence, for all $M$, $\Phi_{\pQs_{\bZ}}(\pi_1^*M)$ is generated by objects of the form 
\[
\Phi_{ \pi^\ast \left( L_{\lambda} V \boxtimes L_\lambda \mathcal Q^\vee \right)}(\pi_1^*M) = \textbf{R}\Gamma(M \otimes L_\lambda \mathcal Q^\vee ) \otimes_{\kk} L_\lambda V
\]
all of which lie in $\lweezer_{d,m}$.  Now, since $\pi_1$ is an affine map, $\Db([Z^{\op{ss}} / \gl])$ is generated by the essential image of $\pi_1^*$.  The result follows.
\end{proof}




\subsection{Truncation operator}

In this section we will see that $\Phi_{\widehat{Q}^+_{Z}}$ has a useful description on $\gl$-representations. Yet before we go deeper into the representation theory we define a truncation operator over our field $ \kk $ of characteristic zero.

\begin{defn}\label{trunc}
	Let $M\in\op{Mod}^{\gl}(\kk[\hom(V,W)])$,  we define the \emph{truncation operator} as follows
	\begin{equation*}
	M_{\geq0}:=\Big(M\otimes\kk[\op{End}(V)]\Big)^{\gl}
	\end{equation*}
\end{defn}

 Recall, further that there is a $\gl\ktimes \gl$-module decomposition 
\begin{equation}\label{theend}
\kk[\op{End}(V)]\cong\bigoplus N_i^\vee\kten N_i,
\end{equation}
where we sum over all irreducible representations of $\gl$ with all positive weights \cite{lie}, these representations are also referred to as \emph{polynomial representations}. Since $\gl$ is linearly reductive over a field of characteristic zero, we may decompose any $\gl$-module $M$ as $M\cong\bigoplus M_i$, where $M_i$ is irreducible and we have the following description of the truncation operator \ref{trunc}:
\begin{rem}\label{rem:Qproj}
Let $M\in\op{Mod}^{\gl}(\kk[\hom(V,W)])$; then decompose $M$ over $\kk$ into irreducibles as 
\begin{equation}\label{theM}
M=\bigoplus_{M_i\text{ irreducible}} M_i.
\end{equation}
Then the truncation operator may be described as follows
\begin{equation*}
M_{\geq0}=\underset{\text{and polynomial}}{\bigoplus_{M_i \text{ irreducible}}} M_i \qedhere
\end{equation*}
\end{rem}

\begin{lem}\label{super:lem0}
For any $M\in\op{Mod}^{\gl}(\kk[\hom(V,W)])$, $M_{\geq0}$ is a $\kk[\hom(V,W)]$-submodule of $M$ and $(\_)_{\geq0}$ is exact.
\end{lem}
\begin{proof}
The exactness of the functor follows since $\gl$ is linearly reductive and thus our operator is just a projection. That $M_{\geq0}$ is a $\kk[\hom(V,W)]$-submodule follows since $\kk[\hom(V,W)]_{\geq0}=\kk[\hom(V,W)]$ since $\kk[\hom(V,W)]$ is a polynomial representation. 
\end{proof}
 To deliver a cleaner picture we define some more notation $Y':=\hom(V,W)$. For the remainder of this subsection we will exploit the commutativity of the following diagram.
\begin{center}\begin{tikzcd}
U^+_{Z}\ar[r,hook,"j"]\ar[d,"q_1|_{U^+_Z}"']&\hom(V,W)\oplus\hom(W',V)\ar[d,"q_1"]\\
U^+_{Y'}\ar[r,hook,"i"']&\hom(V,W)
\end{tikzcd}\end{center}
\begin{lem}\label{super:lem1}
Let $M\in\op{Mod}^{\gl}(\kk[\hom(V,W)])$ then 
\begin{align*}
\Phi_{Q_{Y'}}(M)=M_{\geq0}
\end{align*}
\end{lem}
\begin{proof}
The coaction map defines a morphism 
\begin{align*}
M_{\geq 0} \to (\kk[\op{End}(V)]\otimes M_{\geq 0})^{\gl} \hookrightarrow  (\kk[\op{End}(V)]\otimes M)^{\gl},
\end{align*} which we claim is an isomorphism. Notice that the coaction map lands in\\ $ \kk[\op{End}(V)] \subset \kk[\gl]$ as $ M_{\geq 0} $ is a polynomial representation. To check that this map is an isomorphism, we may base change to $ \overline{\kk} $ (which is faithfully flat over $ \kk $). Hence, assume that $ \kk = \overline{\kk}  $. 

Using Equation~\eqref{theend} and Remark~\ref{rem:Qproj}, we get
\begin{align*}
(\kk[\op{End}(V)]\otimes M)^{\gl} &\cong \bigoplus N_j\otimes(N_j^\vee\otimes M_i ) ^{\gl}\\
&\cong  M_{\geq0}
\end{align*}
where we are considering the left $\gl$ invariant submodule and 
the second line follows from Schur's Lemma.

Finally,  by Lemma \ref{lem:yeesh} we have $Q_{Y'}\cong\kk[A]\otimes\kk[\op{End}(V)]$, and we get
\begin{align*}
(Q_{Y'}\otimes M)^{\gl}&\cong\kk[A]\otimes(\kk[\op{End}(V)]\otimes M)^{\gl}\\
&\cong M_{\geq 0}. \qedhere
\end{align*}
\end{proof}

\begin{lem}\label{super:lem2}
We have an isomorphism
\begin{equation*}
(q_1 \times \op{Id})_* \  {}_s{(Q_{Z})}_p \cong (\op{Id} \times q_1)^* {}_s{(Q_{Y'})}_p 
\end{equation*}
as objects of $\op{Mod}^{{\gl} \times \gl}(Y' \times Z)$.
\end{lem}
\begin{proof}
This follows from the following calculation.
\begin{align*}
_s Q_Z &\cong\kk[A^L,B^R,C] \\
&\cong\kk[A^R,B^R]\otimes_{\kk[A]}\kk[A^L,C] \\
&\cong Z\otimes_{\kk[Y']}Q_{Y'} ,
\end{align*}
where the first isomorphism follows from Lemma \ref{repQ} and in the second line, $k[A]$ acts on the left by going to $A^R$ and on the right by going $CA^L$.
\end{proof}
\begin{cor}\label{super:cor3}
Let $M\in\op{Mod}^{\gl}(k[Y'])$, then 
\begin{equation*}
\Phi_{Q_Z}(q_1^*M) \cong q_1^*\Phi_{Q_{Y'}}(M)
\end{equation*}
\end{cor}
\begin{proof}
This follows from Lemma \ref{super:lem2} which says that it is true at the level of the Fourier-Mukai kernels.
\end{proof}
\begin{lem}\label{thejuice}
For $L_\alpha V\in\lweezer_{d,m}$ we have that 
\begin{align*}
\left({\bf R}i_*{\bf L}i^*L_\alpha V\right)_{\geq0}\cong L_\alpha V
\end{align*}
\end{lem}
\begin{proof}
To see this we will denote the irreducible components as $(\_)_\beta$ where $\beta$ is the highest weight corresponding to the isotypical piece, and by $\beta\geq0$ we denote weights correspond to polynomial representations. 
\begin{align}
\left({\bf R}i_*{\bf L}i^*L_\alpha V\right)_{\geq0}&=\bigoplus_{\beta\geq0} ({\bf R}i_*{\bf L}i^*L_\alpha V)_\beta \notag\\
&\cong\bigoplus_{\beta\geq0}({\bf R}i_*{\bf L}i^*L_\alpha V\otimes L_\beta V^\vee)^{\gl} \notag\\
&\cong \bigoplus_{\beta\geq0}({\bf R}i_*{\bf L}i^*(L_\alpha V\otimes L_\beta V^\vee))^{\gl} \notag\\
&\cong\bigoplus_{\beta\geq0} {\bf R}\Gamma(\op{Gr}(d,W),L_\alpha S\otimes L_\beta S^\vee)\notag \\
&\cong\bigoplus_{\beta\geq0} \Gamma(\op{Gr}(d,W),L_\alpha S\otimes L_\beta S^\vee) \label{line1}\\
&\cong\bigoplus_{\beta\geq0} \Gamma(\op{Hom}(V,W), L_\alpha V\otimes L_\beta V^\vee)^{\gl} \label{line2} \\
&\cong\bigoplus_{\beta\geq0} (\op{Sym}(\op{Hom}(W,V)) \otimes L_\alpha V\otimes L_\beta V^\vee)^{\gl} \notag \\
& \cong\op{Sym}(\op{Hom}(W,V)) \otimes L_\alpha V \label{line3} \\
& \cong \O_{\op{Hom}(V,W)} \otimes L_\alpha V \notag
\end{align}
Equation~\eqref{line1} follows from \cite[Lemma 3.2.a]{Kapp} (this uses the assumption that $L_\alpha V\in\lweezer_{d,m}$ and the fact that the weights of the irreducible summands of $L_\alpha V\otimes L_\beta V^\vee$ are all strictly larger than $-(m-d)$.)  Equation~\eqref{line2} follows as $\op{Gr}(d,W)$ has co-dimension greater than $ 2 $ in the global quotient stack $[\op{Hom}(V,W) / \gl]$. Equation~\eqref{line3} follows from Schur's Lemma and the fact that all representations in $\op{Sym}(\op{Hom}(W,V)) \otimes L_\alpha V$ are polynomial (this uses the fact that $L_\alpha V$ is polynomial).
\end{proof}

\begin{prop}\label{contains2}
If $L_{\alpha}V\in\lweezer_{d,m}$ then 
\begin{equation*}
\Phi_{Q^+_Z}(L_\alpha V)\cong L_\alpha V
\end{equation*}
\end{prop}
\begin{proof}
This result follows from another calculation,
\begin{align*}
\Phi_{Q^+_Z}(L_\alpha V)&\cong \Phi_Q({\bf R}j_*{\bf L}j^*L_\alpha V)\\
&\cong\pi^*\Phi_{Q_{Y'}}({\bf R}i_*{\bf L}i^*L_\alpha V)\\
&\cong \pi^*\Big({\bf R}i_*{\bf L}i^*L_\alpha V\Big)_{\geq0},
\end{align*}
where the second line follows from Corollary \ref{super:cor3} and the last line by Lemma \ref{super:lem1}. Hence our result follows from Lemma \ref{thejuice}.
\end{proof}

\begin{cor} \label{corollary:generators}
$
\op{Im}\Phi_{\widehat{Q}^+_{Z}}=\langle\lweezer_{d,m}\rangle.
$
\end{cor}

\begin{proof}
This is an immediate consequence of Proposition~
\ref{container1} and Lemma \ref{contains2}. 
\end{proof}

Note that we have a similar equality for $\Phi_{\widehat{Q}^-}$.
\begin{cor} \label{cor:imQ-}
	$
	\op{Im}\Phi_{\widehat{Q}^-_{Z}}= \langle \left( \lweezer_{d,m'} \right)^\vee \rangle = \langle \kapreezer_{d,m'} \otimes \op{det}(V^\ast)^{m'-d} \rangle.
	$
\end{cor}

\begin{proof}
 We can switch the roles of $W$ and $W^\prime$ by taking transposes. This is anti-equivariant, i.e., equivariant up to inversion in $\gl$. Consequently, we replace all representations with their duals which gives the first equality. The second is a standard identity. 
\end{proof}

\subsection{The equivalence}

Finally, we combine things to provide Fourier-Mukai equivalences for (twisted) Grassmann flops. As usual, let $\kk$ be an (arbitrary) field of characteristic zero. 

We recall that $P$ is the object obtained by the restriction of $ Q $ to $ Z^+ \times Z^- $.

\begin{thm} \label{theorem:equivalance}
 Assume $\op{dim} W' \geq \op{dim} W$. The wall crossing functor
  \begin{displaymath}
 \Phi_{P} : \Db(Z^+) \to \Db(Z^-)
 \end{displaymath}
 is fully-faithful. If $\op{dim} W' = \op{dim} W$, it is an equivalence. 
\end{thm}

\begin{proof}
 Proposition~\ref{phidecomp} tells us that $\Phi_{\widehat{Q}^+_Z}$ is fully-faithful. Thus, we reduce to checking that $j^\ast_-$ is fully-faithful on the image of $\Phi_{\widehat{Q}^+_Z}$. Also, from Proposition~\ref{phidecompminus}, we know that $j_-^\ast$ is fully-faithful on the image of $\Phi_{\widehat{Q}^-_Z}$. 
 
 From Corollaries ~\ref{corollary:generators} and \ref{cor:imQ-}, we see that 
 \begin{displaymath}
   \operatorname{Im} \Phi_{\widehat{Q}^+_Z} \subseteq \operatorname{Im} \Phi_{\widehat{Q}^-_Z} \otimes \op{det}(V^\ast)^{d-m'}.
 \end{displaymath}
 Since restriction commutes with tensoring with a line bundle, if $j^\ast_-$ is fully-faithful on a full subcategory $\mathcal C$ then it is also on $\mathcal C \otimes \mathcal L$ for any line bundle $\mathcal L$. Now Corollaries~\ref{corollary:generators} and~\ref{cor:imQ-} show $j^\ast_-$ must be fully-faithful on the image of $\Phi_{\widehat{Q}^+_Z}$.  The containment becomes an equality in the case $\op{dim} W' = \op{dim} W$. 
  \end{proof}

\begin{rem} \label{remark:sod}
In Section~\ref{sec:sod}, we use the fully-faithful wall-crossing functors (when, say $ \op{dim} W' > \op{dim} W $) in order to construct semi-orthogonal decompositions for $ \Db(Z^-) $.
\end{rem}

\begin{rem} \label{remark:notDS}
 If $\overline{\kk} = \kk$, once one knows that 
 \begin{displaymath}
  \op{Im} \Phi_{\widehat{Q}^+_Z} = \langle \kapreezer_{d,m} \rangle
 \end{displaymath}
 one can conclude Theorem~\ref{theorem:equivalance} using \cite[Proposition 3.6]{DS}. But, the technology presented here makes for a simple direct proof.
 \end{rem}

\begin{rem} \label{remark:general_k}
 In general, if we have two smooth projective varieties $X$ and $Y$ over $\kk$, then the existence of an equivalence 
 \begin{displaymath}
  \Db(X_{\overline{\kk}}) \cong \Db(Y_{\overline{\kk}})
 \end{displaymath}
 does not guarantee the existence of an equivalence
 \begin{displaymath}
  \Db(X) \cong \Db(Y). 
 \end{displaymath}
 A simple class of counter-examples is Severi-Brauer varieties. 
 
 One needs, at least, a kernel over $\kk$ which base changes to furnish the equivalence to appeal to \cite[Lemma 2.12]{OrlAb}. Without providing a kernel for general $\kk$ for the equivalence in \cite{DS}, the results in loc.cit. cannot be used to deduce equivalences over arbitrary fields of characteristic zero.
\end{rem}

One can go even further. We give the following definition.

\begin{defn}
 We say 
 \begin{center}
 \begin{tikzpicture}[scale=1,level/.style={->=stealth,thick},evel/.style={<->=stealth,thick}]
	\node (a) at (-1.5,1.5) {$Y^+$};
	\node (b) at (1.5,1.5) {$Y^-$};
	\node (c) at (0,0) {$Y_0$};
	\draw[evel,dashed] (a) -- (b) ;
	\draw[level] (a) -- (c) ;
	\draw[level] (b) -- (c) ;
 \end{tikzpicture}
 \end{center}
 is a \textit{twisted Grassmann flop} if the base change to the separable closure of $\kk$
 \begin{center}
 \begin{tikzpicture}[scale=1,level/.style={->=stealth,thick},evel/.style={<->=stealth,thick}]
	\node (a) at (-1.5,1.5) {$Y^+_{\kk^{\op{sep}}}$};
	\node (b) at (1.5,1.5) {$Y^-_{\kk^{\op{sep}}}$};
	\node (c) at (0,0) {$(Y_0)_{\kk^{\op{sep}}}$};
	\draw[evel,dashed] (a) -- (b) ;
	\draw[level] (a) -- (c) ;
	\draw[level] (b) -- (c) ;
 \end{tikzpicture}
 \end{center}
 is isomorphic to a Grassmann flop. 
\end{defn}

\begin{ex}
 Let $A$ be a central simple $\kk$-algebra of degree $n$. For $0 < l < n$, the $l$-th generalized Severi-Brauer variety  $\op{SB}_l(A)$ of $A$ is the variety parameterizing right ideals of dimension $ln$ in $A$. Such a variety is a twisted form of $\op{Gr}(l,n)$, ie 
 \begin{displaymath}
   \op{SB}_l(A)_{\kk^{\op{sep}}} \cong \op{Gr}(l,n)_{\kk^{\op{sep}}}.
 \end{displaymath}
 On $\op{SB}_l(A)$, the tautological vector bundle $\mathcal T$, whose fibers are the ideals, base changes to $\op{Hom}(W,\mathcal S)$. Let $T$ denote the associate geometric vector bundle. The map 
 \begin{displaymath}
   \op{SB}_l(A) \to \op{Spec} \Gamma(T, \mathcal O_T)
 \end{displaymath}
 contracts the zero section and base changes to $Z^+ \to Z_0$. One can then take two copies of $\op{Spec} \Gamma(T, \mathcal O_T)$ and identify them with the involution that base changes to transposition the linear maps. The resulting diagram is a(n honestly) twisted Grassmann flop. 
\end{ex}

We also have equivalences for twisted Grassmann flops in characteristic zero. 

\begin{cor} \label{corollary:twisted}
 Assume $\op{char} \kk = 0$. If we have a twisted Grassmann flop, then there is an equivalence
 \begin{displaymath}
  \Db(Y^+) \to \Db(Y^-).
 \end{displaymath}
\end{cor}

\begin{proof}
 Theorem~\ref{theorem:q_fib_prod} says that the structure sheaf of the fiber product $Y^+_{\kk^{\op{sep}}} \times_{(Y_0)_{\kk^{\op{sep}}}} Y^-_{\kk^{\op{sep}}}$ is a Fourier-Mukai kernel. Applying \cite[Lemma 2.12]{OrlAb} shows that the $Y^+ \times_{Y_0} Y^-$ is also a Fourier-Mukai kernel. 
\end{proof}

\subsection{Semi-orthogonal decompositions} \label{sec:sod}

In this section, we identify the orthogonal to the image of the wall-crossing functor studied in the previous section.  We will assume throughout this section that $ \op{dim} W  = m $ is strictly greater than $ \op{dim} W' = m' $,  and thus obtain a semi-orthogonal decomposition for $ D^b(Z^+) $. (Of course, one can study the case where  $ m' $ is strictly greater than $ m $ using very similar methods.)

Firstly we need to introduce some relevant notation. Consider a vector space $ H $ of dimension $ d - 1 $. Then we have the following morphism
\[
	\overline{\mathcal S} :=	\left[ \Hom(W',H) \oplus \Hom(H,V) / \op{GL}(H) \times \gl \right] \longrightarrow  \left[ \Hom(W',V) / \gl \right]\, ,
\] obtained by composition and then forgetting the $ \op{GL}(H) $-action. We consider the open substack $ \tilde{\mathcal S} \subset \overline{\mathcal S} $ where we restrict to injective maps  $ \Hom(H,V)^{\op{inj}}  \subset  \Hom(H,V) $, and  denote the map from $ \tilde{\mathcal S} $ to $ \left[ \Hom(W',V) / \gl \right] $ as
\[
\tilde{h} : \tilde{\mathcal S} \to  \left[ \Hom(W',V) / \gl \right].
\] By base changing from $  \Hom(W',V)  $ to  $ Z $, we get the morphism
\[
 h : \mathcal S \to \left[ Z/\gl \right],
\] where 
\[
\mathcal S :=	\left[ \Hom(V,W) \oplus \Hom(W',H) \oplus \Hom(H,V)^{\op{inj}}   / \op{GL}(H) \times \gl \right].
\]

\begin{rem}\label{rem:fiber}
 Note that the fiber of the map $ \tilde{h} $ (similarly, $ h $) at a point  $ t $ in $ \left[ \Hom(W',V) / \gl \right] $ is empty unless $ t $ is a map that is not of full-rank, i.e., not surjective. The non-trivial fibers are  the quotients $ \left[\Hom(H,\coker t)^{\op{ss}}/\op{GL}(H)\right] $, where the semi-stable locus is the set of maps of rank $ d - \op{dim( Im} t) - 1  $. We may identify this with  the projective space $ \mathbb P^\vee (\coker{t}) \cong \left[\Hom(\coker t, \kk)^{\op{surj}} / \mathbb{G}_m \right] $.
\end{rem}

Now, we define certain subcategories of $ D(\op{QCoh}^{\gl}(Z)) $, which we will later identify as orthogonals to our window subcategories.

\begin{defn}\label{def:O}
	We define the subcategory  $ \mathsf O_{ d,s} $, where $ d < m' \leq s $, of $ D(\op{QCoh}^{\gl}(Z)) $ inductively on $ s $.
	\begin{itemize}
		\item For $m' \leq j \leq s$, define $\mathsf {US}_{d,j} $  as the one generated by the objects  $ \left( h_*\left( L_{\lambda} H^\vee\right)\right)^\vee  $, 
		\[
		\mathsf {US}_{d,j} := \left\langle \left( h_*\left( L_{\lambda} H^\vee\right)\right)^\vee \right\rangle,
		\] where $ \lambda $ runs over the set of Young diagrams of height $ j +1  - d $ and width $ d -1 $
		\item Let  $ \mathsf O_{ d, m'}  =   \mathsf {US}_{ d, m'} $.  Then for $j \geq m'$, define  $ \mathsf O_{ d, j + 1} $ as the smallest thick triangulated subcategory generated by the objects $ \mathsf O_{ d,j}  \otimes\, \det V$ and  $\mathsf {US}_{d,j}$. 
	\end{itemize} 
\end{defn}

\begin{rem}
The category $\mathsf O_{ d, s}$ is the smallest thick triangulated subcategory generated by the objects  $\mathsf {US}_{ d, i} \otimes \det V^{s-i}$ for $m' \leq i \leq s$.
\end{rem}

The semi-orthogonal decomposition of $ 	\Db(Z^+) $ is given by the following theorem.

	\begin{thm}
	Let $d \leq m' \leq s$.
		There is a semi-orthogonal decomposition
		\[
		\kapreezer_{d,s+1} = \langle   \mathsf O_{d, s},  \kapreezer_{d,m'} \rangle
		\]
				which induces a semi-orthogonal decomposition
		\[
		\Db(Z^+)  = \langle \mathsf O_{ d,m-1}, \Db(Z^-) \rangle
		\]
	\end{thm}
\begin{proof}
By Corollary~\ref{cor:imQ-} and Proposition~\ref{phidecomp}, $\kapreezer_{d,m'}$ is admissible.  Hence, it suffices to check generation and orthogonality.

	Lemma~\ref{lem:gen} shows that $ \kapreezer_{d,j+1}$ 	  is generated by the categories $ \kapreezer_{d,j} $ and $ \mathsf O_{d,j} $. By induction, we see that $  \kapreezer_{d,s} $ is generated by $ \mathsf O_{d,s}, \mathsf O_{d,s-1}, \cdots, \mathsf O_{d,m'}     $ and $ \kapreezer_{d,m'} $.  However, notice  that $ \mathsf O_{d,s} \supset \mathsf O_{d,s-1} \supset \cdots \supset  \mathsf O_{d,m'} $ from Definition~\ref{def:O}, and hence we see that $ \kapreezer_{d,s+1} $ is generated by $ \mathsf O_{d,s} $ and $ \kapreezer_{d,m'} $. The orthogonality of $ \mathsf O_{d,s} $ and $ \kapreezer_{d,m'} $ is the statement of Lemma~\ref{lem:orth}.
	
	In order to get the semi-orthogonal decomposition of $ \Db(Z^+) $, we use the equivalences provided by Corollary~\ref{corollary:generators} and Corollary~\ref{cor:imQ-}.
\end{proof}

\begin{rem}
	As mentioned earlier, the map $ h $ lands in the negative unstable locus of $ Z $, i.e., where the maps $ \Hom(W',V) $ are not surjective. Hence the generating objects of $ \mathsf O_{d,m'}  $ are supported on the negative unstable locus as well.
\end{rem}

We claim that the category $ \kapreezer_{d,m} $ is generated by the categories $ \kapreezer_{d,m'} $ and $  \mathsf O_{d,m'}  $. In order to prove this, we need to recall certain exact sequences discovered by Donovan and Segal in  \cite[Appendix A.2]{DS}.

\begin{prop}[{\cite[Theorem A.7]{DS}}]\label{prop:dsexact}
	Let $\delta$ be a Young diagram of width $  < d $. Then there is an exact sequence of sheaves
	\begin{equation}\label{eq:exact}
	 0 \to L_{\delta^K} V^\vee \otimes \bigwedge^{s_K} W' \to \ldots \to  L_{\delta^1} V^\vee \otimes \bigwedge^{s_1} W' \to  L_{\delta^0} V^\vee \to h_*\left( L_\delta H^\vee\right) \to 0\, ,
	\end{equation}
	where  the $ \delta^k $ and $ s_k $ are defined as follows.
	
	We define a sequence of Young diagrams $ \delta^k $ starting from $ \delta = \delta^0 $  of width $  < d $: 
	\begin{itemize}
		\item $ \delta^1 $ is obtained from $ \delta^0 $ by adding boxes to the first row until it reaches width $ r $. 
		\item $ \delta^k $ is obtained from $ \delta^{k-1} $ by adding boxes to the $ k $-th row until its width is one more than the width of the $ (k-1) $-th row of $ \delta$.
 	\end{itemize}
 Then, $ s_k $ is defined as the difference in the size of the diagrams $ \delta^k $ and $ \delta^0 $. The sequence  terminates when we reach a positive integer $ K $ such that $ s_{K+1} > m'$.
\end{prop}

\begin{proof}
    We refer the reader to \cite[Theorem A.7]{DS}. We note that the restriction on the height of the Young diagram in the statement of Theorem A.7 is unnecessary, as the statement is proved in Section A.3 of \textit{loc.cit.} without any such restrictions.
\end{proof}

\begin{rem}
	We note that the notational difference in the above proposition (where columns and rows have been exchanged) to Theorem A.7 in \cite{DS} can be attributed to the difference in definitions of our Schur functor $ L_\lambda $ to the $ \mathbb{S}^\lambda $ of Donovan and Segal. The contents of both statements are exactly the same.
\end{rem}

In what follows, we use the following standard fact about $\gl$-representations \cite[\S 2 Exercise 18]{wey}.
\begin{prop} 
\label{prop: exer18}
There is a canonical isomorphism
\[
L_{\lambda_1, ..., \lambda_s} V = L_{n- \lambda_s, ..., n- \lambda_1} V^\vee \otimes \det V^s. 
\]
\qed \end{prop}

\begin{lem} \label{lem:gen}
    Assume $d \leq m' \leq s$. 
    The category $ \kapreezer_{d,s+1} $ is generated by the categories $ \kapreezer_{d,s} $ and $  \mathsf O_{d,s} $.
\end{lem}

\begin{proof}
    The idea of the proof is as follows: we can partition  the Young diagrams of 
    $ \kapreezer_{d,s+1} $ by the number of full rows and use the exact sequence \eqref{eq:exact} to work a downward induction on that number.

    Let us set 
    \begin{displaymath}
        \Lambda_t := \left\{ L_\lambda V \left| \lambda_i = d \text{ for } 1 \leq i \leq s-m'+1-t\,, \lambda_{s+1-d} \neq 0 \text{ and } \lambda_i = 0 \text{ for } i > s-d+1 \right. \right\}.
    \end{displaymath}
    So $\Lambda_t$ has $s-m'+1-t$ full rows of length $d$. Furthermore, we set $\Lambda_{-1} := \kapreezer_{d,s}$. 
    For each $\lambda \in 
    \Lambda_t$ we set 
    \begin{displaymath}
        \lambda = (\underbrace{d,\cdots,d}_{s-m'+1-t \op{times}}, \delta).
    \end{displaymath}
   
    The base case of $t=-1$ is clear. Next, we treat the inductive step for $0 \leq t < s-m'+1$. To do so, we tensor the exact sequence~\eqref{eq:exact} for $\delta$ by $ \left( \det V^{\vee }\right)^{s-m'} $. From Lemma~\ref{prop: exer18}, tensoring $ L_\nu V^\vee $ with $ \det V^\vee $ gives $ L_{\nu'} V^\vee $, where $ \nu' $ is obtained by  adding a row of size $ d $ to the diagram $ \nu $.
    Then all the terms appearing in the exact sequence, excluding the last term $ h_*\left( L_\delta H^\vee\right) \otimes \left( \det V^{\vee }\right)^{s-m'}  $, and the first term $ L_{\lambda} V^\vee $, only involve diagrams that appear in $ \Lambda_s $ with $s < t$. Dualizing, we get triangles generating $L_{\lambda} V$ from $ \Lambda_s $ with $s < t$ and $\mathsf O_{d,s} $.

    Finally, we treat the inductive step of $t = s-m'+1$, i.e., no full rows. 
	Choose $ \lambda  $ to be a Young diagram that is a part of the definition of $ \kapreezer_{d,s+1} $ of width $ d-1 $ that is not in $ \kapreezer_{d,s} $. Then, all the terms except for the right-most two terms in the exact sequence~\eqref{eq:exact} (note $\delta = \lambda$ here) belong to the set of vector bundles that we have already generated under the induction hypothesis. Again, by taking the dual of the sequence, we see that we can generate $ L_\lambda V$.
\end{proof}

In order to get a semi-orthogonal decomposition of $ \kapreezer_{d,m'+1} $, we need the following cohomology vanishings.
\begin{lem}\label{lem:orth}
	There are no $ \Hom $s of any homological degree from  $ \kapreezer_{d,s}  $ to $ \mathsf{O}_{d,s} $, i.e.
	\[
		\Hom^\ast(  \kapreezer_{d,s}, \mathsf{O}_{d,s}  ) = 0.
	\]
\end{lem}

\begin{proof}
	Consider a generator $ \left(h_*\left( L_\delta H^\vee\right)\right)^\vee \otimes\, \left(\det V\right)^{\otimes s-i} $ of $ \mathsf{O}_{d,s} $, where $ 0 \leq s-i \leq s+1-d-h(\delta) $ ($ h(\delta) $ denotes the height of $ \delta $), and a generator $ L_\lambda V $ of $ \kapreezer_{d,m'}  $. Then, we want to show that the following $ \Hom $ set vanishes (for all shifts).
	\begin{align}
	&\Hom^\ast \left( L_\lambda V , \left(h_*\left( L_\delta H^\vee\right)\right)^\vee \otimes \left(\det V\right)^{\otimes s-i} \right)  \\
	&= \Hom^\ast \left(  h_*\left( L_\delta H^\vee\right), L_\lambda V^\vee \otimes \left(\det V\right)^{\otimes s-i} \right) \nonumber\\
	&=  \Hom^\ast \left(   L_\delta H^\vee, h^! L_\lambda V^\vee \otimes \left(\det V\right)^{\otimes s-i} \right) \nonumber \\ \label{eq:5.8}
	&= \mathbf{R}\Gamma\left(\mathcal S, L_\delta H \otimes  L_\lambda V^\vee \otimes \left(\det V\right)^{\otimes m'+1-d} \otimes \left(\det V\right)^{\otimes s-i} \otimes \left(\det {H^\vee}\right)^{\otimes m'-d} \right)\\ \label{eq:5.9}
	&= \mathbf{R}\Gamma\left(\mathcal S, L_\delta H \otimes  L_{\lambda'} V  \otimes \left(\det {H^\vee}\right)^{\otimes m'-d} \otimes \det V^{s-i} \right),
	\end{align}	In order to get line~\eqref{eq:5.8}, we use the following expression \cite[Equation (28)]{DS} for the upper shriek functor
	\[
		h^!(-) = h^*(-) \otimes \left(\det V\right)^{\otimes m'+1-d} \otimes \left(\det {H^\vee}\right)^{\otimes m'-d} [ d-m'-1].
	\] In line~\eqref{eq:5.9},  $\lambda'  := (\underbrace{d,\cdots,d}_{s-i+1 \op{times}}, d - \lambda_{m'-d}, ..., d - \lambda_1)$ and we use the identity of Proposition~\ref{prop: exer18}. In particular, note that $ \lambda' $ has at least one row of length $ d $.

	We claim that all the cohomologies in  line~\eqref{eq:5.9} vanish. Notice first that the stack 
	\[
	\left[\Hom(W',H) \oplus \Hom(H,V)^{\op{inj}}/\op{GL}(H) \right]
	\] is the total space of a vector bundle over the projective space $ \mathbb P(V^\vee) = \op{Gr}(d-1,V) $, and (after suppressing the $ \gl $-invariants)  line~\eqref{eq:5.9} becomes
	\begin{align}
	& \mathbf{R}\Gamma\left(\left[\Hom(W',H) \oplus \Hom(H,V)^{\op{inj}}/\op{GL}(H) \right],  L_\delta H \otimes  L_{\lambda'} V  \otimes \left(\det {H^\vee}\right)^{\otimes m'-d}   \otimes \op{Sym} (W^\vee \otimes V)\right) \nonumber \\ 
	\label{eq:coh}
	&= \mathbf{R} \Gamma\left( \op{Gr}(d-1,V), L_\delta H \otimes  L_{\lambda'} V  \otimes \left(\det {H^\vee}\right)^{\otimes m'-d}  \otimes \op{Sym} (W^\vee \otimes V)\otimes \op{Sym}( W' \otimes H^\vee) \right) \\ \label{eq:514}
	&=  \mathbf{R} \Gamma\left( \op{Gr}(d-1,V), L_{\delta'} H^\vee  \otimes \left(\det {H}\right)^{\otimes h(\delta) - m' + d} \otimes  L_{\lambda'} V   \otimes \op{Sym} (W^\vee \otimes V)\otimes \op{Sym}( W' \otimes H^\vee) \right) \\ \label{eq:515}
	&=  \bigoplus_{\mu} \mathbf{R}\Gamma\left( \op{Gr}(d-1,V), L_{\delta'} H^\vee \otimes L_\mu H^\vee   \otimes \left(\det {H}\right)^{\otimes  h(\delta) - m' + d} \otimes  L_{\lambda'} V   \otimes \op{Sym} (W^\vee \otimes V)\otimes L_\mu W' \right) \\ \label{eq:516}
	&= \bigoplus_{\mu,\nu} \mathbf{R}\Gamma\left( \op{Gr}(d-1,V),( L_{\nu} H^\vee)^{\oplus c^\nu_{\delta', \mu} }  \otimes \left(\det {H}\right)^{\otimes  h(\delta) - m' + d} \otimes  L_{\lambda'} V   \otimes \op{Sym} (W^\vee \otimes V)\otimes L_\mu W' \right).
	\end{align} Here,  $ h(\delta) $ denotes the height of the diagram $ \delta $ and  \hbox{$ \delta ' = (d-1-\delta_{h(\delta)}, \cdots, d-1 - \delta_1) $}.   In line~\eqref{eq:514} we use Proposition~\ref{prop: exer18}.    In line~\eqref{eq:515}, we use the decomposition
	\[
	\op{Sym} (A \otimes B) = \bigoplus_\mu L_\mu(A) \otimes L_\mu(B),
	\] where we sum over all Young diagrams $ \mu $. To get the last line~\eqref{eq:516}, we use the Littlewood-Richardson rule, and  $ c^{\nu}_{\delta',\mu} $ denotes the Littlewood-Richardson coefficients.
	
		Now, we can use the tautological exact sequence on $ \op{Gr}(d-1,V) $ to get the following identification
	\[
	H^\vee = \mathcal U^\vee, \qquad \det H = \mathcal U^{\perp} \otimes \det V,
	\] where $ \mathcal U^\vee $ and $ \mathcal U^{\perp} $ are the vector bundles appearing in the tautological exact sequence
	\[
	0 \to \mathcal U^{\perp} \to V^\vee \otimes \mathcal O_{\op{Gr}(d-1,V)} \to \mathcal U^\vee \to 0.
	\] This allows us to use the Borel-Weil-Bott theorem \cite{Dem} to compute the cohomologies appearing in equation~\eqref{eq:516} (see \cite{kuz} for a review of this method) as follows. 
	
	First, we assume that $ h(\delta) - m' + d \leq 0 $. Then the expression in 	equation~\eqref{eq:516} becomes
	\[
	= \bigoplus_{\mu,\nu} \Gamma \left( \op{Gr}(d-1,V),( L_{\nu'} H^\vee)^{\oplus c^\nu_{\delta', \mu} }  \otimes  L_{\lambda'} V   \otimes \op{Sym} (W^\vee \otimes V)\otimes L_\mu W' \right).
	\]
	where $\nu'$ is obtained from $\nu$ by adding $m' - h(\delta) - d$ rows of length $d-1$ to the top of $\nu$.  The only non-trivial vector bundle component in the above expression is $L_{\nu'} H^\vee$ which has no higher cohomology since it is a Schur functor applied to $H^\vee = \mathcal U^\vee$.  Furthermore, global sections of $L_{\nu'} H^\vee = L_{\nu'} \mathcal U^\vee$ returns nothing more than $L_{\nu'} V^\vee$.

Hence, the above expression is reduced to
	\begin{align*}
	 \bigoplus_{\mu,\nu} ( L_{\nu'} V^\vee)^{\oplus c^\nu_{\delta', \mu} } \otimes  L_{\lambda'} V   \otimes \op{Sym} (W^\vee \otimes V)\otimes L_\mu W',
	\end{align*} where the sum is over all $ \nu $ of width $ < d $ (since $\dim H = d-1$).  Now simply notice that $\nu'$ always has width $ < d$ but any irreducible representation in  $L_{\lambda'} V   \otimes \op{Sym} (W^\vee \otimes V)$ is built from adding to $\lambda'$ which has a row length $d$.  Therefore, no representations can cancel and this expression vanishes upon taking $\gl$-invariants.

	Now, we consider the case  when  $   h(\delta) - m' + d > 0 $. Applying Borel-Weil-Bott, see e.g. \cite[Theorem 3.1, Corollary 3.4]{kuz} , gives us the following.
	\begin{align}
	&\bigoplus_{\mu,\nu} \mathbf{R}\Gamma\left( \op{Gr}(d-1,V),( L_{\nu}  \mathcal U^\vee)^{\oplus c^\nu_{\delta', \mu} }  \otimes \left(\mathcal U^{\perp} \otimes \det V\right)^{\otimes  h(\delta) - m' + d} \otimes  L_{\lambda'} V   \otimes \op{Sym} (W^\vee \otimes V)\otimes L_\mu W' \right)\nonumber \\ \label{eq:517}
	&= \bigoplus_{\mu,\nu} ( L_{\sigma\bullet \alpha }  V^\vee  [l(\sigma)])^{\oplus c^\nu_{\delta', \mu} }   \otimes \left( \det V\right)^{\otimes  h(\delta) - m' + d} \otimes  L_{\lambda'} V   \otimes \op{Sym} (W^\vee \otimes V)\otimes L_\mu W' .
    \end{align} 
    Here  we define  the  diagram $ \alpha $ (of width $ d $) by adding a column of height $ h(\delta) - m' + d  $ in the $ d $-th position (the right-most one) to $ \nu $. We define an element $ \sigma $  of the symmetric group $ S_{d} $ (uniquely) as follows. To the diagram $ \alpha $, we first `add' a diagram $ \rho = (d,d-1,\cdots, 1) $ vertically, i.e., we add $ d $ boxes to the first column, $ d-1 $ boxes to the second column and so on. Then, we pick (the unique) symmetric group element $ \sigma $  that permutes the columns in order to make them non-increasing to give a diagram, say $ \sigma.\alpha $. Finally, we subtract $ \rho $ from $ \sigma.\alpha $ vertically, i.e, we remove $ d $ boxes from the first column, $ d-1 $ boxes from the second column and so on. This gives a diagram that we denote by $ \sigma \bullet \alpha $. (For the convenience of the reader, we provide a short graphical illustration of the procedure described above to get the diagram $ \sigma \bullet \alpha $ from the diagram $ \nu $, immediately after the end of this proof.)
	
	The key observation for us is that the diagram $\sigma \bullet \alpha  $ has at most $ h(\delta) - m' + d $ full rows. This is because the $ d $-th column of $ \alpha $ has height $ h(\delta) - m' + d  $, and the operation $ \sigma \bullet \alpha $ cannot increase the height of this column. Thus, we see that every term $ ( L_{\sigma\bullet \alpha }   V^\vee) \otimes \left( \det V\right)^{\otimes  h(\delta) - m' + d} $ appearing in equation~\eqref{eq:517} can be rewritten as $ (L_\gamma V^\vee) \otimes \left(\det V\right)^{p}  $ where $ \gamma $ has no full rows  and $ p \geq 0 $. On the other hand, as before $ L_{\lambda'} V $ has a full row and hence the $ \gl $-invariants of equation~\eqref{eq:517} vanish.
\end{proof}

\begin{ex}

We illustrate the procedure in the above proof to get the diagram $ \sigma \bullet \alpha $ from the diagram $ \nu $  (using the  notation introduced in the proof). Assume $ d = 3 $, and that $ h(\delta) - m' +d $ = 3 and  choose a diagram $ \nu = (2,1) $. Then, we have the following,

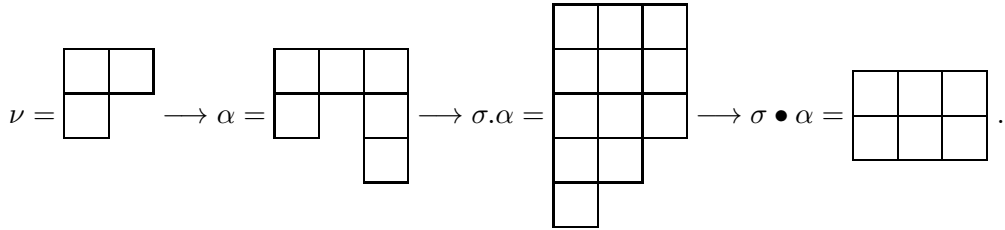
\[
\begin{aligned}
\nu = \begin{ytableau}
\tikznode{a1}{~} &   \\
\tikznode{a1}{~} & \none \\
\none & \none   \\
\end{ytableau}
\longrightarrow
\alpha = 
\begin{ytableau}
\tikznode{a1}{~}   & \tikznode{a2}{~}  &  \\
\tikznode{a1}{~}   & \none  &   \\
\none   & \none  &   \\
\end{ytableau}
\longrightarrow
\sigma . \alpha  = 
\begin{ytableau}
\tikznode{a1}{~}   & \tikznode{a2}{~}  &  \\
\tikznode{a1}{~}   & \none  &   \\
\tikznode{a3}{~}   &   &   \\
\tikznode{a1}{} &  & \none \\
\tikznode{a1}{} & \none & \none 
\end{ytableau} \longrightarrow
\sigma \bullet \alpha  = 
\begin{ytableau}
\tikznode{a1}{~}   & \tikznode{a2}{~}  &  \\
\tikznode{a1}{~}   &  &   
\end{ytableau} 
\end{aligned} \,\,.
\]

\end{ex}

\section*{Acknowledgements}
This work was supported by a grant from the Simons Foundation (708132, MRB). M.\ Ballard and R.\ Vandermolen were partially supported by NSF DMS-1501813. D.\ Favero and N.K. Chidambaram were partially supported by NSERC RGPIN 04596 and CRC TIER2 229953. M. Ballard thanks W. Donovan for stimulating discussions on Grassmann flops. D.\ Favero is greatly appreciative to B. Kim for useful discussions related to this work.  We thank the anonymous referee for various comments and suggestions that resulted in various improvements, both in presentation and mathematical content.


\begin{thebibliography}{BLVdB16}
	
	\bibitem[BDF17]{BDF}
	M.~R. {Ballard}, C.~{Diemer}, and D.~{Favero}.
	\newblock {Kernels from Compactifications}.
	\newblock {\em ArXiv e-prints}, October 2017.
	
	\bibitem[BFK14]{BFKe}
	Matthew Ballard, David Favero, and Ludmil Katzarkov.
	\newblock A category of kernels for equivariant factorizations and its
	implications for hodge theory.
	\newblock {\em Publications math{\'e}matiques de l'IH{\'E}S}, 120(1):1--111,
	Nov 2014.
	
	\bibitem[BFK19]{BFK}
	Matthew Ballard, David Favero, and Ludmil Katzarkov.
	\newblock Variation of geometric invariant theory quotients and derived
	categories.
	\newblock {\em J. Reine Angew. Math.}, 746:235--303, 2019.
	
	\bibitem[BLVdB16]{BLV}
	Ragnar-Olaf Buchweitz, Graham~J. Leuschke, and Michel Van~den Bergh.
	\newblock Non-commutative desingularization of determinantal varieties, {II}:
	{A}rbitrary minors.
	\newblock {\em Int. Math. Res. Not. IMRN}, (9):2748--2812, 2016.
	
	\bibitem[BO95]{BO}
	A.~{Bondal} and D.~{Orlov}.
	\newblock {Semiorthogonal decomposition for algebraic varieties}.
	\newblock In {\em eprint arXiv:alg-geom/9506012}, June 1995.
	
	\bibitem[Dem76]{Dem}
	Michel Demazure.
	\newblock A very simple proof of {B}ott's theorem.
	\newblock {\em Invent. Math.}, 33(3):271--272, 1976.
	
	\bibitem[{Dri}13]{drin}
	V.~{Drinfeld}.
	\newblock {On algebraic spaces with an action of $\mathbb{G}_m$}.
	\newblock {\em ArXiv e-prints}, August 2013.
	
	\bibitem[DS14]{DS}
	Will Donovan and Ed~Segal.
	\newblock Window shifts, flop equivalences and {G}rassmannian twists.
	\newblock {\em Compos. Math.}, 150(6):942--978, 2014.
	
	\bibitem[HL15]{HL}
	Daniel Halpern-Leistner.
	\newblock The derived category of a {GIT} quotient.
	\newblock {\em J. Amer. Math. Soc.}, 28(3):871--912, 2015.
	
	\bibitem[{Kap}88]{Kapp}
	M.~M. {Kapranov}.
	\newblock {On the derived categories of coherent sheaves on some homogeneous
		spaces.}
	\newblock {\em Inventiones Mathematicae}, 92:479, 1988.
	
	\bibitem[Kaw02]{Kaw}
	Yujiro Kawamata.
	\newblock D -equivalence and k -equivalence.
	\newblock {\em J. Differential Geom.}, 61(1):147--171, 05 2002.
	
	\bibitem[KP96]{kraft}
	H.~Kraft and C.~Procesi.
	\newblock {\em A Primer of Classical Invariant Theory: Preliminary Version}.
	\newblock {\'e}diteur inconnu, 1996.
	
	\bibitem[Kra10]{kra}
	H.~Krause.
	\newblock Localization theory for triangulated categories.
	\newblock {\em London Mathematical Society Lecture Note Series}, 2010.
	
	\bibitem[Kuz08]{kuz}
	Alexander Kuznetsov.
	\newblock Exceptional collections for {G}rassmannians of isotropic lines.
	\newblock {\em Proc. Lond. Math. Soc. (3)}, 97(1):155--182, 2008.
	
	\bibitem[Mum65]{GIT}
	David Mumford.
	\newblock {\em Geometric invariant theory}.
	\newblock Ergebnisse der Mathematik und ihrer Grenzgebiete, Neue Folge, Band
	34. Springer-Verlag, Berlin-New York, 1965.
	
	\bibitem[Orl02]{OrlAb}
	D.~O. Orlov.
	\newblock Derived categories of coherent sheaves on abelian varieties and
	equivalences between them.
	\newblock {\em Izv. Ross. Akad. Nauk Ser. Mat.}, 66(3):131--158, 2002.
	
	\bibitem[Pro07]{lie}
	Claudio Procesi.
	\newblock {\em Lie Groups: An Approach through Invariants and Representations}.
	\newblock Universitext. Springer Science+Buisness Media LLC, 2007.
	
	\bibitem[SVdB17]{SV}
	\v{S}pela Spenko and Michel Van~den Bergh.
	\newblock Non-commutative resolutions of quotient singularities for reductive
	groups.
	\newblock {\em Invent. Math.}, 210(1):3--67, 2017.
	
	\bibitem[Wey46]{wey46}
	Hermann Weyl.
	\newblock {\em The Classical Groups, Their Invariants and Representations, by
		Hermann Weyl,...[with Supplement.].}
	\newblock Princeton University Press, 1946.
	
	\bibitem[Wey03]{wey}
	Jerzy Weyman.
	\newblock {\em Cohomology of Vector Bundles and Syzygies}.
	\newblock Cambridge Tracts in Mathematics. Cambridge University Press, 2003.
	
\end{thebibliography}
\end{document}